\documentclass{elsart}
\usepackage{natbib,amsmath,amssymb,amscd}
\usepackage[dvips]{graphicx}

\newcommand{\romannum}{\renewcommand{\labelenumi}{\textnormal{(\roman{enumi}).}}}

\newcommand\vphi{\varphi}
\newcommand\vrho{\varrho}
\newcommand\vtheta{\vartheta}

\newcommand\eps{\varepsilon}
\newcommand\nats{\mathbb{N}}

\newcommand\reals{\mathbb{R}}
\newcommand\cplxs{\mathbb{C}}

\newcommand\ball{\mathbb{B}}

\newcommand\vvoid{\varnothing}
\newcommand\sle{\leqslant}
\newcommand\sge{\geqslant}
\newcommand\limk{\lim\nolimits}
\newcommand\equiva{\Leftrightarrow}

\newcommand\clos[1]{\overline{{#1}}}
\newcommand\var{\llcorner\kern-.3em\lrcorner}

\newcommand\mathtxt[1]{\quad\text{{#1}}\quad}
\newcommand\smathtxt[1]{\ \text{{#1}}\ }

\newcommand{\nd}{\mathtxt{and}}

\newcommand\fs{for some }

\newcommand\Fa{For all }
\newcommand\fa{for all }
\newcommand\mathfa[1][{}]{\quad\text{\fa{#1} }}
\newcommand\smathfa[1][{}]{\ \text{\fa{#1} }}

\newcommand\sth{so that }
\newcommand\scth{such that }

\DeclareMathOperator\dom{dom}
\DeclareMathOperator\epi{epi}

\DeclareMathOperator\ind{ind}
\DeclareMathOperator\id{id}
\DeclareMathOperator\Index{Index}
\DeclareMathOperator\Gr{Gr}

\makeatletter
\def\Set@Scallop[#1]#2#3{{#1}\Parens{#2}{#3}}
\newcommand\DeclareScalableOperator[2]{%
  \expandafter\def\csname#1\endcsname{\@ifnextchar[{{#2}\Set@Scallop}{{#2}\Set@Scallop[{}]}}
}
\newcommand\Size[7][1]{
                                 \ifx#20%
                                        \def\r@l{}\def\r@m{}\def\r@r{}%
                                 \else%
                                    \ifx#21%
                                           \def\r@l{\bigl}\def\r@r{\bigr}\def\r@m{\bigm}%
                                    \else%
                                           \ifx#22%
                                                 \def\r@l{\Bigl}\def\r@r{\Bigr}\def\r@m{\Bigm}%
                                            \else%
                                                 \ifx#23%
                                                        \def\r@l{\biggl}\def\r@r{\biggr}\def\r@m{\biggm}%
                                                  \else
                                                        \ifx#24%
                                                              \def\r@l{\Biggl}\def\r@r{\Biggr}\def\r@m{\Biggm}%
                                                        \fi%
                                                  \fi%
                                            \fi%
                                      \fi%
                                 \fi%
                                 \ifx#10%
                                       \def\r@m{}%
                                 \fi%
                                 \r@l#3{#4}\r@m#5{#6}\r@r#7%
}%
\newcommand\Set[3]{
                                 \Size{#1}{\{}{#2}{|}{#3}{\}}
}%
\newcommand\Rscp[3]{
                                 \Size[0]{#1}{(}{#2}{:}{#3}{)}%
}%
\newcommand\Parens[2]{
  \Size[0]{#1}{(}{#2}{}{}{)}
}
\newcommand\Bracks[2]{
  \Size[0]{#1}{[}{#2}{}{}{]}
}

\newcommand\Lopint[3]{
  \Size[0]{#1}{]}{#2}{,}{#3}{]}%
}
\newcommand\Ropint[3]{
  \Size[0]{#1}{[}{#2}{,}{#3}{[}%
}
\newcommand\Norm[2]{
  \Size[0]{#1}{\lVert}{#2}{}{}{\rVert}
}
\newcommand\Abs[2]{
  \Size[0]{#1}{\lvert}{#2}{}{}{\rvert}
}
\newcommand\Span[2]{
  \Size[0]{#1}{\langle}{#2}{}{}{\rangle}
}

\DeclareScalableOperator{Lip}{\mathrm{Lip}} 
\DeclareScalableOperator{Mea}{\mathfrak{M}} 
\DeclareScalableOperator{Ct}{\mathcal{C}} 
\DeclareScalableOperator{Cc}{\mathcal{C}_c} 
\DeclareScalableOperator{Lp}{\mathbf{L}} 
\DeclareScalableOperator{Hp}{\mathbf{H}} 
\DeclareScalableOperator{Cst}{\mathrm{C}^*} 
\DeclareScalableOperator{Cred}{\mathrm{C}^*_r} 
\DeclareScalableOperator{Cpt}{\mathbb{K}} 
\DeclareScalableOperator{Lct}{\mathcal{L}} 
\DeclareScalableOperator{Mult}{\mathrm{M}} 
\DeclareScalableOperator{Hom}{\mathrm{Hom}} 
\DeclareScalableOperator{Diff}{\mathrm{Diff}} 

\DeclareScalableOperator{Co}{\mathrm{co}}
\DeclareScalableOperator{faces}{\mathcal F}
\DeclareScalableOperator{closed}{\mathbb F}

\newif\if@smallmat
\newif\if@none
\newif\if@paren
\newif\if@brack
\newif\if@brace
\newif\if@vline
\newenvironment{Matrix}[2][1]
                                 {\ifx#20%
                                        \@smallmattrue%
                                  \else%
                                         \@smallmatfalse
                                  \fi%
                                  \ifx#11%
                                         \@nonefalse\@parentrue\@brackfalse\@bracefalse\@vlinefalse%
                                  \else%
                                       \ifx#12%
                                            \@nonefalse\@parenfalse\@bracktrue\@bracefalse\@vlinefalse%
                                        \else%
                                            \ifx#13%
                                                 \@nonefalse\@parenfalse\@brackfalse\@bracetrue\@vlinefalse%
                                            \else%
                                                 \ifx#14%
                                                       \@nonefalse\@parenfalse\@brackfalse\@bracefalse\@vlinetrue
                                                 \else%
                                                       \ifx#15%
                                                             \@nonefalse\@parenfalse\@brackfalse\@bracefalse\@vlinefalse%
                                                       \else%
                                                             \@nonetrue\@parenfalse\@brackfalse\@bracefalse\@vlinefalse%
                                                       \fi%
                                                 \fi%
                                            \fi%
                                        \fi%
                                   \fi%
                                   \if@smallmat%
                                        \if@none%
                                             \begin{smallmatrix}%
                                        \else%
                                            \if@paren%
                                                  \bigl(\begin{smallmatrix}%
                                            \else%
                                                  \if@brack%
                                                          \bigl[\begin{smallmatrix}%
                                                  \else%
                                                          \if@brace%
                                                               \bigl\{\begin{smallmatrix}%
                                                          \else%
                                                               \if@vline%
                                                                    \bigl\lvert\begin{smallmatrix}%
                                                                \else%
                                                                    \bigl\lVert\begin{smallmatrix}%
                                                                \fi%
                                                          \fi%
                                                  \fi%
                                            \fi%
                                        \fi%
                                   \else%
                                        \if@none%
                                             \begin{matrix}%
                                        \else%
                                            \if@paren%
                                                  \begin{pmatrix}%
                                            \else%
                                                  \if@brack%
                                                          \begin{bmatrix}%
                                                  \else%
                                                          \if@brace%
                                                               \begin{Bmatrix}%
                                                          \else%
                                                               \if@vline%
                                                                    \begin{vmatrix}%
                                                                \else%
                                                                    \begin{Vmatrix}%
                                                                \fi%
                                                          \fi%
                                                  \fi%
                                            \fi%
                                        \fi%
                                   \fi}%
                                  {\if@smallmat%
                                        \if@none%
                                             \end{smallmatrix}%
                                        \else%
                                            \if@paren%
                                                  \end{smallmatrix}\bigr)%
                                            \else%
                                                  \if@brack%
                                                          \end{smallmatrix}\bigr]%
                                                  \else%
                                                          \if@brace%
                                                               \end{smallmatrix}\bigr\}%
                                                          \else%
                                                               \if@vline%
                                                                    \end{smallmatrix}\bigr\rvert%
                                                                \else%
                                                                    \end{smallmatrix}\bigr\rVert%
                                                                \fi%
                                                          \fi%
                                                  \fi%
                                            \fi%
                                         \fi%
                                   \else%
                                        \if@none%
                                             \end{matrix}%
                                        \else%
                                            \if@paren%
                                                  \end{pmatrix}%
                                            \else%
                                                  \if@brack%
                                                          \end{bmatrix}%
                                                  \else%
                                                          \if@brace%
                                                               \end{Bmatrix}%
                                                          \else%
                                                               \if@vline%
                                                                    \end{vmatrix}%
                                                                \else%
                                                                    \end{Vmatrix}%
                                                                \fi%
                                                          \fi%
                                                  \fi%
                                            \fi%
                                        \fi%
                                   \fi}%

\def\thmskip{\vskip3\p@ plus 2\p@ minus 2\p@}
\newdimen\l@belindent\l@belindent\z@
\def\m@kel@bel#1{\hspace{\m@rgin}\hspace{\labelwidth}\hspace{\labelsep}\hspace{\l@belindent}#1}
\def\redefine@makelabel{\let\m@rgin=\leftmargin\setlength{\leftmargin}{\z@}\def\makelabel{\m@kel@bel}}
\newcounter{the@rem}
\def\thethe@rem{(\roman{the@rem})}
\newenvironment{multipr@@f}[1][\empty]%
{\par\penalty-500\ignorespaces\thmskip%
  \begin{list}{\textbf{PROOF of~\thethe@rem{.}}\ignorespaces}%
  {\usecounter{the@rem}\listparindent=\parindent\itemsep=0.5ex\topsep=3\p@\l@belindent\z@\redefine@makelabel}}%
{\quad\hfill\quad\ \qed%
  \end{list}\ignorespacesafterend\par\penalty-600}
\newenvironment{Partsproof}[1][\empty]{\ignorespaces\begin{multipr@@f}[#1]}{\end{multipr@@f}\ignorespacesafterend}
\newcommand\incrementcounter{\addtocounter{the@rem}1}
\makeatother

\begin{document}

\begin{frontmatter}

\title{Spectrum and Analytical Indices of the C$^*$-Algebra of Wiener--Hopf Operators}
\author{Alexander Alldridge\corauthref{corauth}},
\ead{alldridg@math.upb.de}
\corauth[corauth]{Corresponding author. Address: Institut f\"ur Mathematik, Universit\"at Paderborn, Germany.}
\author{Troels Roussau Johansen\thanksref{thks}}
\ead{johansen@math.upb.de}
\thanks[thks]{Part of the work was carried out while the author was supported
by IHP `Harmonic Analysis and Related Problems', HPRN-CT-2001-00273. He is currently supported by a DFG post.doc grant under
the International Research \& Training Group ``Geometry and
Analysis of Symmetries'', {\texttt{http://irtg.upb.de}}}
\address{Universit\"at Paderborn, Institut f\"ur Mathematik, Warburger Strasse 100, D--33098 Paderborn, Germany}
\date{\today}
\maketitle
\thispagestyle{empty}

\begin{abstract}
We study multivariate generalisations of the classical Wiener--Hopf algebra, which is the C$^*$-algebra generated by the Wiener--Hopf operators, given by the convolutions restricted to convex cones. By the work of Muhly and Renault, this C$^*$-algebra is known to be isomorphic to the reduced C$^*$-algebra of a certain restricted action groupoid, given by the action of Euclidean space on a certain compactification. Using groupoid methods, we construct composition series for the Wiener--Hopf C$^*$-algebra by a detailed study of this compactification. We compute the spectrum, and express homomorphisms in $K$-theory induced by the symbol maps which arise by the subquotients of the composition series in analytical terms. Namely, these symbols maps turn out to be given by an analytical family index of a continuous family of Fredholm operators. In a subsequent paper, we also obtain a topological expression of these indices. 
\end{abstract}

\end{frontmatter}

\section{Introduction}

\subsection{The Classical Wiener--Hopf Equation}

 The classical Wiener--Hopf equation is of the form $(1+W_f)u=v\,$, where
  \[
  W_fu(x)=\int_0^\infty f(x-y)u(y)\,dy\mathfa f\in\Lp[^1]0\reals\,,\,u\in\Lp[^2]0{0,\infty}\,,\,x\in\Ropint00\infty
  \ .
  \]
  The bounded operator $W_f$ is called the Wiener--Hopf operator of symbol $f\,$. The operator $W_f$ is conjugate, via the Euclidean
  Fourier transform, to the Toeplitz operator $T_{\hat f}$ defined on the Hardy space of the upper half plane, and thus has
  connections to both complex and harmonic analysis. Moreover, its multi-variate generalisations (see below) play a role in applications to wave propagation, e.g.~in the
  presence of diffraction by an impenetrable wedge-shaped obstruction.

  The one-variable Wiener--Hopf equation is well understood, by the following classical theorem.

\begin{thm}[\cite{gohberg-krein}]
Let $\mathcal W(0,\infty)$ be the C$^*$-algebra generated by the operators $W_f\,$, $f\in\Lp[^1]0\reals\,$.
\begin{enumerate}\romannum
\item The following sequence is exact
\[
\begin{CD}%
0@>>>{\Lct1{\protect{\Lp[^2]0{0,\infty}}}}@>>>{\mathcal W(0,\infty)}@>\sigma>>{\Ct[_0]0\reals}@>>>0%
\end{CD}
\]
where $\sigma$ is the Wiener--Hopf representation, defined by $\sigma(W_f)=\Hat f\,$.
\item The operator $1+W_f$ is Fredholm if and only if $1+\Hat f$ is everywhere non-vanishing on $\reals^+\,$.
\item In this case, $\Index(1+W_f)$ is the negative winding number of $(1+\Hat f)(e^{i\theta})\,$.
\end{enumerate}
\end{thm}

\subsection{Multivariate Generalisation}

  It is quite straightforward to generalise the above setting to several variables. Indeed, let $X$ be a finite-dimensional real
  vector space endowed with some Euclidean inner product $\Rscp0\var\var\,$, and let $\Omega\subset X$ be a closed, pointed and solid
  convex cone. I.e., $\Omega$ contains no line, and has non-void interior. Consider Lebesgue measure on $X$ to define $\Lp[^1]0X$ and
  its restriction to $\Omega$ to define $\Lp[^2]0\Omega\,$.

  Then, $W_f$ is defined by
  \[
  W_fu(x)=\int_\Omega f(x-y)u(y)\,dy\mathfa f\in\Lp[^1]0X\,,\,u\in\Lp[^2]0\Omega\,,\,x\in\Omega\ .
  \]
  Moreover, let $\mathcal W(\Omega)$ be the \emph{Wiener--Hopf algebra}, the C$^*$-subalgebra of all bounded operators on
  $\Lp[^2]0\Omega$ generated by the collection of the $W_f\,$, $f\in\Lp[^1]0\Omega\,$.

  The programme we propose to study then is the following:
  {\def\theenumi{\arabic{enumi}}
  \begin{enumerate}
  \item Determine a composition series of $\mathcal W(\Omega)$ and compute its subquotients.
  \item Find Fredholmness criteria for Wiener--Hopf operators. (This is just a reformulation of the computation of subquotients.)
  \item Give an index formula which expresses their Fredholm index in topological terms.
  \end{enumerate}}

  \noindent
  These problems have been addressed from different angles in a quite extensive literature. Pioneering work was done in the series of
  papers \cite{douglas-coburn-wh1,douglas-coburn-wh2}, 
  \cite{douglas-coburn-schaeffer-singer,douglas-coburn-singer}.
  Together with the work of \cite{douglas-howe}, this culminated in
  the solution of problems (1)-(3) for the example of the (discrete) quarter plane. \cite{berger-coburn-u2} were the first to address
  the structure of the Hardy--Toeplitz algebra (equivalent to the Wiener--Hopf algebra for symmetric tube type domains) for a symmetric domains of
  rank $2\,$, the $2\times2$ matrix ball (the rank $1$ case having been essentially solved by Venugopalkrishna). This led to the paper
  of \cite{berger-coburn-koranyi} which treats the case of all Lorentz cones (also corresponding to rank $2$ symmetric domains, the Lie
  balls).

  Major advances were made in \cite{upmeier-solvable,upmeier-wh,upmeier-fredholm} who solved (1) for the Hardy--Toeplitz algebras of all bounded
  symmetric domains (which properly include the Wiener--Hopf algebras for symmetric cones). Moreover, he developed an index theory,
  proving index formulae for the all Wiener--Hopf operators associated to symmetric cones, and thus solving problem (3) for this class of
  cones. A basic tool in his approach is the Cayley transform, which allows for the transferral to the situation of bounded symmetric
  domains.

  Another approach was taken by \cite{dynin-wh1}, who uses an inductive procedure, based on the local decomposition of the cone $\Omega$ into
  a product relative to a fixed exposed face, to construct a composition series as in (1). This presumes
  a certain tameness of the cone $\Omega\,$, which he calls `complete tangibility'. Due to the weakness of this assumption,
  a large class of cones, including polyhedral, almost smooth and homogeneous cones, are subsumed.

  \smallskip\noindent
  The point of view we will adopt in this note is due to \cite{muhly-renault}. They describe a general procedure to produce a
  (locally compact, measured) groupoid whose groupoid C$^*$-algebra is the Wiener--Hopf algebra, and compute composition series (1) for the opposite extremes of polyhedral
  and symmetric cones. Their construction is based on the specification of a convenient compactification of $\Omega$ (in fact, of $X$). \cite{nica-remarks-cone}
  has given a uniform construction of this Wiener--Hopf compactification for \emph{all} pointed and solid cones. The main problem is to prove that the corresponding groupoid 
  always has a Haar system. From the more general perspective of ordered homogeneous spaces, in which $X$ is replaced by a locally compact group
  and $\Omega$ by a submonoid satisfying certain assumptions, \cite{hn-wiener1} have extended Nica's results, at the same time giving a
  convenient alternative description of the Wiener--Hopf compactification.

  As yet, none of the problems (1)-(3) have been solved in full generality. In fact, there is not even an index theorem for the
  polyhedral case. We show how the groupoid perspective allows for a unified treatment, for a large class of cones satisfying some global regularity assumption
  which arises in a natural fashion.

  \medskip\noindent
  In this paper, we obtain a composition series of the Wiener--Hopf algebra, in the following manner. The Wiener--Hopf algebra is isomorphic to the reduced groupoid
  C$^*$-algebra $\Cred0{\mathcal W_\Omega}$ of the groupoid $\mathcal W_\Omega\,$, defined as the restriction $(\clos X\rtimes X)|\clos\Omega$ of the action groupoid 
  given by the action of the vector space $X$ on $\clos X\,$, the order compactification of  $X$ (see below), to the closure $\clos\Omega$ of $\Omega$ in $\clos X\,$.

  Order the dimensions of convex faces of the dual cone $\Omega^*$ increasingly by
  \[
  \{0=n_0<n_1<\dotsm<n_d=n\}=\Set1{\dim F}{F\subset\Omega^*\text{ face}}\ .
  \]
  Let $P_j$ be the set of $n_{d-j}$-dimensional faces of $\Omega^*\,$, and assume
  that it is compact \fa $j\,$, in the space of all closed subsets of $X\,$, endowed with the Fell topology. (This class of cones properly contains the polyhedral and
  symmetric cases, where the $P_j$ are, respectively, finite sets and certain compact homogeneous spaces including, in particular, all spheres.) Then there is a
  surjection from $\clos\Omega=\mathcal W_\Omega^{(0)}$ onto the set of all faces of $\Omega^*$ which is continuous when restricted to the inverse image $Y_j$ of
  $P_j\,$.

  The sets $Y_j$ are closed and invariant, and $U_j=\bigcup_{i=0}^{j-1}Y_i$ are open and invariant. Thus, we obtain ideals
  $I_j=\Cred0{\mathcal W_{\Omega}|U_j}$ of the Wiener--Hopf C$^*$-algebra $\Cred0{\mathcal W_\Omega}\,$, and extensions
  \[
  \begin{CD}%
    0@>>>\Cred0{\mathcal W_\Omega|Y_{j-1}}@>>>I_{j+1}/I_{j-1}=\Cred0{\mathcal W_\Omega|(U_{j+1}\setminus U_{j-1})}@>>>%
    \Cred0{\mathcal W_\Omega|Y_j}@>>>0\ .
	\end{CD}
  \]
  Moreover, we have Morita equivalences $\mathcal W_\Omega|Y_j\sim\Sigma_j$ where $\Sigma_j=\mathcal W_\Omega|P_j$ is the `co-tautological' topological vector
  bundle over the space $P_j$ whose fibre at the face $F$ is the orthogonal complement $F^\perp\,$. We prove the following theorem.

  \begin{thm}
    The Wiener--Hopf algebra admits an ascending filtration by ideals $(I_j)_{j=0,\dotsc,d}$ whose subquotients are stably isomorphic to $\Ct[_0]0{\Sigma_j}\,$, 
    and in particular, is solvable of length $d$ in the sense of \cite{dynin-invsing}. 
  \end{thm}

  \noindent
  The spectrum can be computed in terms of a suitable gluing of the bundles $\Sigma_j\,$. As a particular case, we obtain the classical Wiener--Hopf extension
  (i.e., $X=\reals$ and $\Omega=\reals_{\sge0}$). 

  Moreover, the above extensions give rise to index maps $\partial_j:K^1_c(\Sigma_j)\to K^0_c(\Sigma_{j-1})\,$.
  In this paper, we give an analytical expression of the $\partial_j\,$, as follows.

  \begin{thm}
    The quotient $I_{j+1}/I_j$ is a field $\mathcal K_j$ of elementary C$^*$-algebras over $\Sigma_j\,$. If a class $f\in K^1_c(\Sigma_j)$ is represented by a element invertible
    modulo matrices over $I_j\,$, then its image $\sigma_j$ in the matrices over $\Mult0{\mathcal K_j}$ is a Fredholm family, and
    \[
    \partial_j(f)=\Index_{\Sigma_{j-1}}\sigma_j(f)
    \]
    is the analytical family index.
  \end{thm}

  \noindent
  In the second part of our work, \cite{alldridge-johansen-wh2}, we also give a topological expression of $\partial_j\,$, which we now proceed to explain.
  Assume that the cone $\Omega$ has a facially compact and locally smooth dual cone (compare section 6 of \cite{alldridge-johansen-wh2}). Consider the compact
  space $\mathcal P_j$ consisting of all pairs $(E,F)\in P_{j-1}\times P_j$ \scth $E\supset F\,$. It has (not necessarily surjective) projections
  \[
  \begin{CD}P_{j-1}@<\xi<<\mathcal P_j@>\eta>>P_j\end{CD}
  \]
  The map  $\xi:\mathcal P_j\to P_{j-1}$ turns $\mathcal P_j$ into a fibrewise $\mathcal C^1$ manifold over the compact base $\xi(\mathcal P_j)\,$. Moreover,
  $\eta^*\Sigma_j$ is the trivial line bundle over $\xi^*\Sigma_{j-1}\oplus T\mathcal P_j$ if $T\mathcal P_j$ denotes the fibrewise tangent bundle. Then we have the
  following theorem

  \begin{thm}
    The $KK$-theory element representing $\partial_j$ is given by
    \[
    \partial_j\otimes\zeta^*=\eta^*\otimes y\otimes\tau_j\mathtxt{in}KK^1(\Cred0{\Sigma_j},\Cred0{\Sigma_{j-1}|\xi(\mathcal P_j)})\ ,
    \]
    where $y\in KK^1(\cplxs,S)$ represents the classical Wiener--Hopf extension, $\eta^*$ is associated to the projection $\eta^*\Sigma_j\to\Sigma_j\,$, and $\zeta^*$ is
    associated to the inclusion $\Sigma_{j-1}|\xi(\mathcal P_j)\subset\Sigma_{j-1}\,$.

    Here, $\tau_j\in KK(\Cred0{\xi^*\Sigma_{j-1}\oplus T\mathcal P_j},\Cred0{\Sigma_{j-1}|\xi(\mathcal P_j)})$ represents the topological Atiyah--Singer family index for
    $\xi^*\Sigma_{j-1}\oplus T\mathcal P_j\,$, considered as a vector bundle over $\Sigma_{j-1}|\xi(\mathcal P_j)\,$.
  \end{thm}

\section{The Wiener--Hopf C$^*$-Algebra}

\subsection{The Wiener--Hopf Groupoid}
Let $X$ be a finite-dimensional Euclidean vector space, $\Omega\subset X$ a closed convex cone which we assume to be \emph{pointed} ($-\Omega\cap\Omega=0$) and 
\emph{solid} ($\Omega-\Omega=X$).

  In order to construct a groupoid which conveniently describes the C$^*$-algebra of Wiener--Hopf operators, we recall the order compactification of the Euclidean space
  $X\,$. Here, we follow \cite{hn-wiener1}. (The compactification was first described in \cite{nica-remarks-cone}, in a quite different manner.)

  Consider the set $\closed0X$ of closed subsets of $X\,$. The \emph{topology of Painl\'ev\'e--Kuratowski convergence} is a complete, compact
  and separable metric topology for which the convergent sequences $(A_k)\subset\closed0X$ are those for which
  $\varlimsup_kA_k=\varliminf_kA_k\,$. Here, the \emph{Painl\'ev\'e--Kuratowski limes inferior} resp.~\emph{superior} are
  \begin{align*}
    \varliminf\nolimits_kA_k&=\bigcap_{\eps>0}\bigcup_{k\in\nats}\bigcap_{\ell\sge k}\Parens1{A_\ell}_\eps
    =\Set1{a=\limk_ka_k\in X}{a_k\in A_k}\\
    \intertext{and}
    \varlimsup\nolimits_kA_k&=\bigcap_{\eps>0}\bigcap_{k\in\nats}\bigcup_{\ell\sge k}\Parens1{A_\ell}_\eps
    =\Set1{a=\limk_ka_{\alpha(k)}\in X}{a_{\alpha(k)}\in A_{\alpha(k)}}\ ,
  \end{align*}
  where $A_\eps=\Set1{x\in X}{\inf\nolimits_{a\in A}\Norm0{x-a}<\eps}\,$.

  Since $X$ is a locally compact, separable metric space, the Painl\'ev\'e--Kuratowski convergence topology coincides with the Fell topology, and the Vietoris topology on the
  one-point compactification $X^+$ \cite[prop.~I.1.54, th.~I.1.55]{hu-papageorgiou-vol1}, i.e.~the topology induced by the Hausdorff metric of $X^+\,$.
  This is the fashion in which the topology on $\closed0X$ was introduced in \cite{hn-wiener1}.

  We inject $X$ into $\closed0X$ by $\eta:X\to\closed0X:x\mapsto x-\Omega\,$. The map $\eta$ is a homeomorphism onto its image, and $\eta(X)$ is open in its closure
  \cite[lem.~II.8, th.~II.11]{hn-wiener1}. Take $\clos X$ to be the closure of $X$ in $\closed0X\,$, and similarly, denote the closure of $\Omega$ by $\clos\Omega\,$.

  The elements of $\clos\Omega$ are non-void, and for any $a\in\clos X\,$, $a\neq\vvoid\,$, there exists $x\in X$ such that $a+x\in\clos\Omega\,$; i.e., $\clos\Omega$ 
  intersects every orbit of $\clos X\,$, save one, under the action of $X\,$. Consequently, $\clos\Omega$ completely determines $\clos X\,$.
\bigskip

 Define a right action of $X$ on $\closed0X$ by $A.x=A+x\,$. Clearly, this action leaves $X$ invariant. Hence, it also leaves $\clos X$ invariant, and we may form the
  transformation groupoid $\clos X\rtimes X\,$. We define $\mathcal W_\Omega=(\clos X\rtimes X)|\clos\Omega\,$, the \emph{Wiener--Hopf groupoid}, as its restriction. Recall
  \[
  r(\omega,x)=\omega\nd s(\omega,x)=\omega+x\mathfa(\omega,x)\in\mathcal W_\Omega\ .
  \]
  The locally compact groupoid $\mathcal W_\Omega$ is topologically amenable. Indeed, $\clos X$ is an amenable $X$-space, since the group $X$ is amenable
  \cite[cor.~2.2.10]{anantharaman-delaroche-renault}. Moreover, $\mathcal W_\Omega$ is topologically equivalent to the restriction of $\clos X\rtimes X$ to the non-void
  elements of $X$ \cite[ex.~2.7 and p.~16]{muhly-renault-williams-equivalences}, and amenability is preserved under topological equivalence
  \cite[th.~2.2.17]{anantharaman-delaroche-renault} and restriction to open invariant subsets.

  A Haar system of $\clos X\rtimes X$ is given by $\lambda^A=\delta_A\otimes\lambda$, $\lambda$ denoting the Lebesgue measure on $X\,$. That this Haar system restricts
  to $\mathcal W_\Omega$ is a non-trivial matter related to the regularity of the compactification $\clos X\,$. It was proved in \cite{muhly-renault} for polyhedral and
  symmetric cones, in \cite[prop.~1.3]{nica-remarks-cone} for the present setup and subsequently in \cite[lem.~III.4]{hn-wiener1}, by a different method, for the more general
  setting of ordered homogeneous spaces.

\begin{thm}[Muhly--Renault, Nica, Hilgert--Neeb]
  The C$^*$-algebra $\Cst0{\mathcal W_\Omega}$ of the locally compact group\-oid $\mathcal W_\Omega$ is isomorphic to the C$^*$-algebra of Wiener--Hopf operators.
\end{thm}

  Moreover, the above authors also established that the ideal of compact operators on $\Lp[^2]0\Omega$ is naturally contained in $\Cst0{\mathcal W_\Omega}\,$. In order
  to describe $\Cst0{\mathcal W_\Omega}$ in greater detail, we embark on a closer study of the compactification $\clos\Omega\,$.

\subsection{The Wiener--Hopf Compactification}
  As a motivating example, consider the quarter plane $\Omega=\Ropint001^2\,\subset\reals^2=X\,$. This cone is self-dual and simplicial.
  Identifying a point $x\in\Omega$ with the set $x-\Omega\,$, we see that limits of sequences $x_k$ can contribute to
  $\clos\Omega\setminus\Omega$ in two distinct fashions. Either, one of the components of $x_k$ remains bounded; in this case, the
  limit point will be an affine half space not completely containing $\Omega\,$. Or, both components tend to infinity; in which case,
  the limit shall be the entire space $X\,$. This is illustrated below.
  \begin{center}
    \includegraphics[bb=148 514 531 669,scale=0.6]{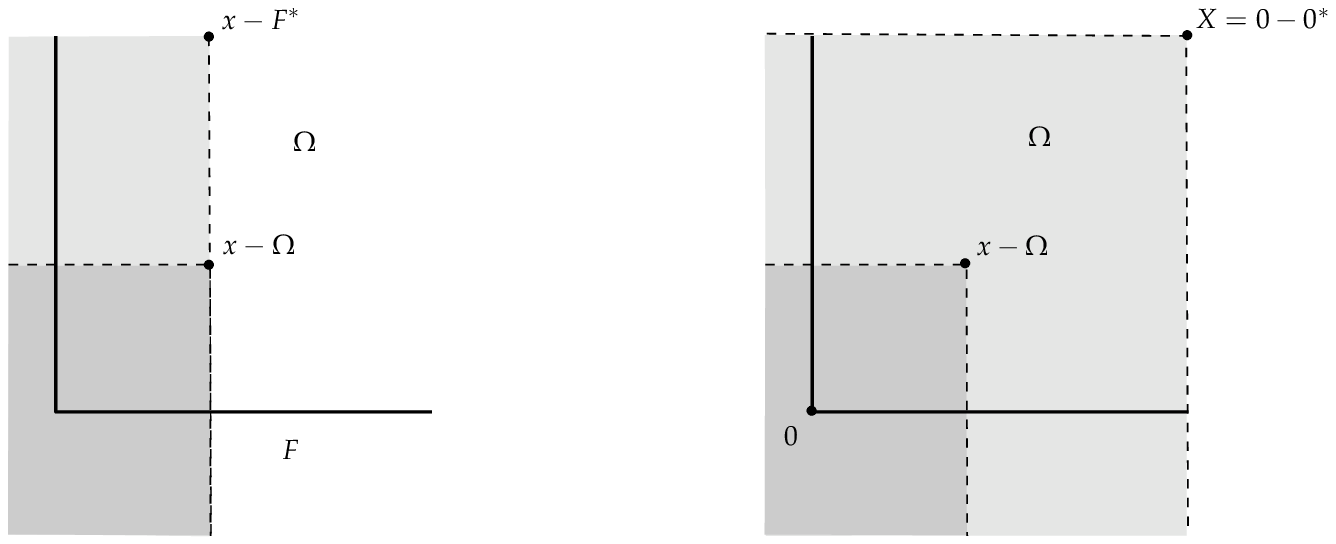}
  \end{center}

  This example suggests that $\clos\Omega$ is a local fibre bundle over the spaces of faces. More precisely, denote by
  \[
  P=\Set1{\vvoid\neq F\subset\Omega^*}{F\text{ convex face}}
  \]
  and \fa $A\subset X\,$, let
  \[
  \Span0F=A-A\ ,\ A^*=\Set1{x\in X}{x\in\Rscp0xA\sge0}\nd A^\circledast=\Span0A\cap A^*
  \]
  denote the linear span, the dual cone, and the relative dual cone, respectively. Then $\clos\Omega=\Set1{x-F^*}{x\in F^\circledast\,,\,F\in P}\,$, at least in the example.

  In general, it seems to be a non-trivial matter to give a complete description of $\clos\Omega\,$. \cite[prop.~4.6.2]{nica-remarks-cone} has proved that at least the
  inclusion $\supset$ in the above equality always holds. For the reader's convenience, we give a streamlined proof, whilst claiming no originality on our part.

\begin{thm}[Nica]\label{compact-bundle}
  We have the inclusion $\Set1{x-F^*}{x\in F^\circledast\,,\,F\in P}\subset\clos\Omega\,$, where $P$ is the face lattice of $\Omega^*\,$. Moreover, $x$ and $F$ are uniquely
  determined by $x-F^*\,$.
\end{thm}

\par\noindent
  We first note the following lemmata.

\begin{lem}\label{face-projection}
  Let $F\subset\Omega$ be an exposed face, i.e.~the intersection of $\Omega$ with a supporting hyperplane. Denote by $F^\perp$ the orthogonal complement, and by
  $\check F=F^\perp\cap\Omega^*$ the dual face. Then
  \[
  \check F^*=\clos{\Omega-F}=\Span0F\oplus\clos{p_{F^\perp}(\Omega)}\nd\clos{p_{F^\perp}(\Omega)}\cap\Span0{\check F}=\check F^\circledast
  \]
  In particular, $\clos{p_{\check F}(\Omega)}=\check F^\circledast\,$. Here, for $A\subset X\,$, $p_A:X\to\Span0A\sqsubset X$ denotes orthogonal projection onto
  the \emph{linear span} of $A\,$.
\end{lem}

\begin{pf}
  The equality $\check F^*=\clos{\Omega-F}$ follows from \cite[prop.~I.1.9]{hilgert-hofmann-lawson}. Since the intersection of the closed convex cones $\Span0F$ and
  $\clos{p_{F^\perp}(\Omega)}$ is $0\,$, their sum is closed \cite[I.2.32]{hilgert-hofmann-lawson}. For the first assertion, it remains to
  prove that $\Omega-F=\Span0F+p_{F^\perp}(\Omega)\,$. But this follows from
  \[
  p_{F^\perp}(\Omega)\subset\Omega+\Span0F=\Omega-F\nd\Omega\subset p_{F^\perp}(\Omega)+\Span0F\ .
  \]

  As to the second assertion, $\clos{p_{F^\perp}(\Omega)}\subset\check F^*\,$, so $\clos{p_{F^\perp}(\Omega)}\cap\Span0{\check F}\subset\check F^\circledast\,$. Conversely,
  for any $x\in\check F^\circledast\subset\check F^*\,$, by the first assertion, $x+f\in\clos{p_{F^\perp}(\Omega)}\subset F^\perp$ \fs $f\in\Span0F\,$. But then
  $f\in F^\perp-x=F^\perp\,$, since $x\in\Span0{\check F}\subset F^\perp\,$. This implies $f=0\,$, so $x\in\clos{p_{F^\perp}(\Omega)}\,$.

  Finally, from the second assertion, we have $\check F^\circledast\subset p_{\check F}\Parens1{\clos{p_{F^\perp}(\Omega)}}\subset\clos{p_{\check F}(\Omega)}\,$, since
  of course $p_{\check F}p_{F^\perp}=p_{\check F}$ and $\Span0{\check F^\circledast}=\Span0{\check F}\,$. Here, the latter follows from the fact that 
  $\check F^\circledast$ is the dual of the pointed cone $\check F\,$, taken in the vector space $\Span0{\check F}\,$.  Conversely,
  \[
  \Rscp1{p_{\check F}(\Omega)}{\check F}=\Rscp1{p_{F^\perp}(\Omega)}{\check F}\subset\Rscp1{\check F^*}{\check F}\subset\reals_{\sge0}\ ,
  \]
  so $\clos{p_{\check F}(\Omega)}\subset\check F^*\cap\Span0{\check F}=\check F^\circledast\,$. This proves the lemma.
\end{pf}

\begin{lem}\label{ray-limit}
  Let $x\in\Omega\,$, and denote by $F=\Omega\cap\Parens1{\Omega^*\cap x^\perp}^\perp$ the exposed face of $\Omega$  generated by $x\,$.
  Then $\limk_{\lambda\to\infty}\Parens1{\lambda\cdot x-\Omega}=\clos{F-\Omega}=-F^*$ in $\closed0X\,$.
\end{lem}

\begin{pf}
  Let $\lambda>0$ and $\omega\in\Omega\,$. Then
  \[
  \lambda\cdot x-\omega=\limk_{\mu\to\infty}\Parens1{(\lambda+\mu)\cdot x-(\omega+\mu\cdot x)}\in\varliminf\nolimits_{\mu\to\infty}
  \Parens1{\mu\cdot x-\Omega}\ ,
  \]
  and consequently
  \[
  \varlimsup\nolimits_{\lambda\to\infty}\Parens1{\lambda\cdot x-\Omega}\subset\clos{\bigcup_{\lambda>0}\Parens1{\lambda\cdot x-\Omega}}
  \subset\varliminf\nolimits_{\lambda\to\infty}\Parens1{\lambda\cdot x-\Omega}\ .
  \]
  This shows that
  \[
  \lim_{\lambda\to\infty}\Parens1{\lambda\cdot x-\Omega}=\clos{\reals_{\sge0}\cdot x-\Omega}=\clos{F-\Omega}=-\check F^*\ ,
  \]
  where we have used Lemma \ref{face-projection}. This gives our contention.
\end{pf}
{\def\Elproofname{\protect{PROOF of Theorem \ref{compact-bundle}.}}
\begin{pf}
  Let $F\in P\,$, $x\in F^\circledast\,$. By \cite[prop.~1.3~(iv)]{hn-liesemigrp}, there are faces $F=F_0\subset\dotsm F_m=\Omega^*\,$,  $F_j$ exposed in $F_{j+1}\,$. Suffices
  to show that all $u-F_j^*\,$, $u\in F_j^\circledast\,$, are contained in the closure of $\Set0{v-F_j^*}{v\in F_{j+1}^\circledast}\,$. By induction, we may assume $m=1\,$,
  i.e.~that $F$ is exposed.

  Taking $y$ in the relative interior of $F^\circ\,$, we obtain
  $-F^*=\clos{\check F-\Omega}=\limk_{\lambda\to\infty}\Parens1{\lambda\cdot y-\Omega}$ by Lemma \ref{ray-limit}, $F$ being exposed.
  Since $x\in F^\circledast=\clos{p_F(\Omega)}$ by Lemma \ref{face-projection}, there exist $y_k\in\Omega\,$, such that $x=\limk_kp_F(y_k)\,$.
  But then $x-F^*=\limk_k\Parens1{y_k+\lambda_k\cdot y-\Omega}$ for any $\lambda_k\to\infty\,$. Thus, $x-F^*\in\clos\Omega\,$.

  As to uniqueness, recall that for $\vvoid\neq A\subset X\,$, the support functional of $A$ is defined by
  \[
  \sigma_A(x)=\sigma_{\Co0A}(x)=\sup\nolimits_{y\in A}\Rscp0xy\in[0,\infty]\ .
  \]
  The equality of sets $u-E^*=v-F^*$ entails the equality of support functionals, so
  \[
  E=\dom\sigma_{u-E^*}=\dom\sigma_{v-F^*}=F\ .
  \]
  In particular, $u,v\in\Span0E\,$. Moreover,
  \[
  \Rscp0ue=\sigma_{u-E^*}(e)=\sigma_{v-E^*}(e)=\Rscp0ve\mathfa e\in E\ .
  \]
  This proves $u=v\,$, and hence, the assertion.
\end{pf}}

\par\noindent
  The natural question to ask is when the above theorem gives a complete description of $\clos\Omega\,$. \cite[prop.~6.1]{nica-remarks-cone} shows that this is the case
  for the rather restricted class of \emph{tame} cones. (Polyhedral cones, Lorentz cones, and cones of dimension $\sle3$ are tame, but symmetric cones coming from irreducible
  Jordan algebras of rank $\sge2$ are not.) On the other hand, he gives an example of a four-dimensional cone where equality fails, cf.~\cite[ex.~5.3.5]{nica-remarks-cone}.

  Nica's work suggests that the description of $\clos\Omega$ is related to the compactness of $P\,$, considered as a subset of $\closed0X\,$. More precisely, we have the
  following theorem, the proof of which the remainder of this section is devoted to.

\begin{thm}\label{compact-bundle2}
Order the face dimensions increasingly,
  \[
  \{n_0<n_1<\dotsm<n_d\}=\Set1{\dim F}{F\in P}\ ,
  \]
  where we set $\dim A=\dim\Span0A$ for $A\subset X\,$. Let $P_j=\{\dim=n_{d-j}\}\,$,
  \[
  Y_j=\Set1{x-F^*}{x\in F^\circledast\,,\,F\in P_j}\ ,\nd U_j=\bigcup_{i=0}^{j-1}Y_i\ .
  \]
  Moreover, define projections $\pi:U_{d+1}\to P$ and $\lambda:U_{d+1}\to X$ by
  \[
  \pi(x-F^*)=F\nd\lambda(x-F^*)=x\ .
  \]
  Then the following holds.]
\begin{enumerate}
\romannum
\item If $\clos\Omega=U_{d+1}\,$, then $P$ is compact.
\item If $P$ is compact, then $(\pi,\lambda)|Y_j$ is continuous \fa $0\sle j\sle d\,$.
\item If $P_j$ is compact \fa $0\sle j\sle d\,$, then $\clos\Omega=U_{d+1}\,$.
\end{enumerate}
\end{thm}

\begin{rem}
  A least if $P$ is modular, $F\mapsto\dim F$ is the rank function of the lattice $P\,$. Thus, the condition that $P_j$ be compact \fa $j$ simply means that $P$ has a 
  continuous rank function. Moreover, if all $P_j$ are compact, $\clos\Omega$ is locally a fibre bundle. We shall prove this more precise statement  below.
\end{rem}

\par\noindent
The following observation is fundamental, albeit elementary.

\begin{prop}\label{lower-dim-closed}
Let $C_k,C\in\closed0X\,$.
\begin{enumerate}
\romannum
\item We have $\dim\varliminf_kC_k\sle\liminf_k\dim C_k\,$.
\item Assume that $C_k,C$ are convex cones, and that $C_k\to C\,$. Then $\Span0C\subset\varliminf_k\Span0{C_k}\,$. Moreover,
  $\Span0{C_k}\to\Span0C$ if and only if $\dim C=\limk_k\dim C_k\,$.
\end{enumerate}
\end{prop}

\begin{Partsproof}
\item Let $q=\dim C$ and choose $x_1,\dotsc,x_q\in C\,$, linearly independent. For $C_k$ close to $C\,$,
  there exist $x_{ik}\in C_k\,$, $x_{ik}$ close to $x_i$ for $i=1,\dotsc,q\,$. But then the matrices $(x_{1k},\dotsc,x_{qk})$ and
  $(x_1,\dotsc,x_q)$ are close. Since the rank function on $\Hom0{\reals^q,X}$ is l.s.c., the rank of these matrices is $q$ for $k$
  sufficiently large. Thus, $q$ is eventually a lower bound for $\dim C_k\,$, and therefore $q\sle m\,$, proving the lemma.

\item Since $\Span0C=C-C\,$, $\Span0{C_k}=C_k-C_k\,$, the inclusion $\Span0C\subset\varliminf_k\Span0{C_k}$ is trivial. Moreover,
  $\varliminf_k\Span0{C_k}$ is manifestly a linear subspace of $X\,$. By the first part, its dimension is $\sle\dim C=\limk_k\dim C_k$
  whenever this limit exists, so in that case, $\Span0C=\varliminf_k\Span0{C_k}\,$.

  Now, let $x_{\alpha(k)}\in C_{\alpha(k)}$ converge to $x\in X\,$, where $\alpha(\nats)$ is cofinal in $\nats\,$. Since $C_{\alpha(k)}\to C\,$, the above shows that
  \[
  x=\limk_kx_{\alpha(k)}\in\varliminf\nolimits_k\Span0{C_{\alpha(k)}}=\Span0C\ .
  \]
  So we have proved that $\varlimsup_k\Span0{C_k}\subset\Span0C\,$.

  As for the converse, let $\Span0C=\limk_k\Span0{C_k}\,$. Since $X$ is locally compact, the Attouch--Wets topology on $\closed0X\setminus\{\vvoid\}$
  coincides with the topology of Painl\'ev\'e--Kuratowski convergence, by \cite[ex.~5.1.10]{beer-topologies-book}. Now, the dual cones $\Span0C^*=C^\perp\,$,
  $\Span0{C_k}^*=C_k^\perp\,$, so by continuity of polarity, \cite[cor.~7.2.12]{beer-topologies-book}, $C^\perp=\limk_kC_k^\perp\,$. The first part implies
  $\dim C\sle\liminf\nolimits_k\dim C_k$ and
  \[
  \dim C=n-\dim C^\perp\sge n-\liminf\nolimits_k\dim C_k^\perp=\limsup\nolimits_k\dim C_k\ ,
  \]
  so $\dim C=\limk_k\dim C_k\,$.
\end{Partsproof}

\begin{cor}\label{proj-cont}
Let $C_k,C\subset X$ be convex cones, \scth $C_k\to C\,$.
\begin{enumerate}
\romannum
\item If $\pi_{C_k}\,$, $\pi_C$ denote the metric projections, cf.~\cite{zarantonello-convex-spectral}, then
  \[
  \pi_{C_k}(x_k)\to\pi_C(x)\mathtxt{whenever} x_k\to x\,,\,x_k,x\in X\ .
  \]
\item Let $\dim C_k=\dim C\,$. Then $p_{C_k}\to p_C$ and $p_{C_k^\perp}\to p_{C^\perp}\,$.
\end{enumerate}
\end{cor}

\begin{Partsproof}
\item Let $y_k=\pi_{C_k}(x_k)\,$. Then $\Norm0{y_k}\sle\Norm0{x_k}\,$, so $y_k$ is bounded, and we may assume that it converges to some $y\in X\,$. Then
  $y\in\varliminf_kC_k=C\,$. Let $u\in C\,$. There exist $u_k\in C_k\,$, $u_k\to u\,$. Thus
  \[
  \Norm0{x-u}=\limk_k\Norm0{x_k-u_k}\sge\Norm0{y_k-u_k}=\Norm0{y-u}\ ,
  \]
  and it follows that $y=\pi_C(x)\,$.
\item By Proposition \ref{lower-dim-closed}, $\Span0{C_k}\to\Span0C\,$. Thus $p_{C_k}\to p_C$ follows from the first part, because $p_C=\pi_{\Span0C}\,$, and we have already noted $C_k^\perp\to C^\perp$ above.
\end{Partsproof}

\par\noindent
  The following lemma constitutes the main step in the theorem's proof. For its proof, recall the following notions: Given a proper l.s.c.~function
  $\vphi:X\to\,\Lopint0{-\infty}\infty\,$, its \emph{epigraph}
  \[
  \epi\vphi=\Set1{(x,y)\in X\times\reals}{y\sge\vphi(x)}
  \]
  is a closed non-void subset of $X\times\reals\,$. Given $(\vphi_k),\vphi$ proper l.s.c., $(\vphi_k)$ is said to \emph{epi-converge}
  to $\vphi$ if $\epi\vphi_k\to\epi\vphi$ in the Painl\'ev\'e--Kuratowski sense (w.r.t.~the box metric on $X\times\reals$).

\begin{lem}\label{projection-subgraph-closed}
Let $F_k\in\faces0{\Omega^*}\,$, $x_k\in F_k^\circledast\,$, $-\Omega\subset C\subset X$ closed and convex, and
  $E\subset\Omega^*$ be closed. Assume that $x_k-F_k^*\to C$ and $F_k\to F\,$.
\begin{enumerate}
\romannum
\item We have $\dom\sigma_C\subset F\,$.
\item If $(x_k)$ is bounded, then $x=\limk_kx_k\in F^\circledast$ exists, $\dom\sigma_C=F\,$,  and $C=x-F^*\,$.
\item If $(x_k)$ is unbounded and $m=\dim F=\limk_k\dim F_k\,$, then $\dim\dom\sigma_C<m\,$. In fact, there exist $E_k,E\in P\,$,
  $E_k\subset F_k\,$, $E\subset F\,$, \scth
  \[
  \dim E_k\,,\,\dim E<m\ ,\quad E_k\to E\nd x_k-E_k^*\to C\ .
  \]
\end{enumerate}
\end{lem}

\begin{Partsproof}
\item For the support functionals, $\vphi_k=\sigma_{x_k-F_k^*}\to\sigma_C=\vphi$ in the sense of epi-convergence, cf.~\cite[cor.~3.4.5]{nica-remarks-cone}. Whenever we have $\vphi(y)<\infty\,$, by \cite[lem.~6.2]{nica-remarks-cone},
  there exist $y_k\in\Omega^*$ such that $y_k\to y$ and $\vphi_k(y_k)\to\vphi(y)\,$. Since $\vphi(y)<\infty\,$, we may assume
  $y_k\in\dom\vphi_k=F_k\,$. Therefore, $y=\limk_ky_k\in\limk_kF_k=F\,$.
\item Let $(x_k)$ be bounded and assume there is some $y\in F$ such that $\vphi(y)=\infty\,$. Then
  by \cite[lem.~6.2]{nica-remarks-cone} $\liminf_k\vphi_k(y_k)\sge\vphi(y)=\infty$ for any $y_k\in F_k\,$, $y_k\to y\,$.
  In particular,
  \[
  \Rscp0{y_{\alpha(k)}}{x_{\alpha(k)}}=\vphi_{\alpha(k)}(y_{\alpha(k)})\to\infty
  \]
  for some subsequence $\alpha\,$. Seeing that $(x_k)$ and $(y_k)$ are bounded, this is a contradiction. Thus, $\dom\vphi=F\,$. Then
  there exists a unique $x\in F^\circledast$ such that $C=x-F^*\,$, by \cite[lem.~6.1]{nica-remarks-cone}. Let $z=\limk_kx_{\alpha(k)}$ be
  any accumulation point of $(x_k)\,$, and take $u\in X\,$. We may write $u=v-w$ for $v,w\in\Omega^*\,$. By
  \cite[lem.~6.2]{nica-remarks-cone}, there exist $v_k,w_k\,$, such that $v_k\to v\,$, $w_k\to w\,$, and
  \begin{align*}
    \Rscp0zu&=\limk_k\Rscp0{x_{\alpha(k)}}{v_{\alpha(k)}-w_{\alpha(k)}}\\
    &=\limk_k\vphi_k(v_k)-\limk_k\vphi_k(w_k)=\vphi(v)-\vphi(w)=\Rscp0xu\ .
  \end{align*}
  Hence, $x=z\,$, and thus $\limk_kx_k=x\,$.
\item Now, consider the case that $(x_k)$ is unbounded and that $\lim\dim F_k=m\,$. Define $y_k=\Norm0{x_k}^{-1}\cdot x_k\,$.
  Passing to a subsequence, we may assume $y=\limk_ky_k$ exists, and $\dim F_k=m$ \fa $k\,$. We have
  $y\in F^*\,$, by continuity of polarity, \cite[cor.~7.2.12]{beer-topologies-book}, and $\Span0F=\limk_k\Span0{F_k}$ by Proposition \ref{lower-dim-closed} (ii). Consequently,
  $y\in F^\circledast\,$.

  The exposed face $E$ of $F^*$ generated by $y$ satisfies
  \[
  -\check E^*=\clos{E-F^*}=\clos{\reals_{\sge0}\cdot y-F^*}\ ,
  \]
  so in order to prove $C-\check E^*=C\,$, it suffices to prove $C+\reals_{\sge0}\cdot y\subset C\,$. Let $\lambda\sge0$ and set
  $\lambda_k=\Norm0{x_k}^{-1}\,$. For any $f_k\in F_k$ such that $x_k-f_k$ converges,
  \[
  x_k-(1-\lambda\lambda_k)\cdot f_k=\lambda\cdot y_k+(1-\lambda\lambda_k)(x_k-f_k)
  \to\lambda\cdot y+\limk_k(x_k-f_k)\ ,
  \]
  and hence $\lambda\cdot y+\limk_k(x_k-f_k)\in\varliminf_k\Parens1{x_k-F_k}=C\,$. Thus, $C-\check E^*=C\,$.

  There exist $y_k\in F_k^\circledast\,$, $y_k\to y\,$. Let $E_k=F_k\cap y_k^\perp\in P\,$. Clearly, $\varlimsup_kE_k\subset F\cap y^\perp\,$.
  Let $f\in F\cap y^\perp\,$. There exist $f_k\in F_k\,$, $f_k\to f\,$. We can write $f_k=e_k+u_k$ with uniquely determined $e_k\in E_k$ and
  $u_k\in-E_k^*\cap\Span0{F_k}\,$. By Corollary \ref{proj-cont} (i), $e_k$ converges to the projection of $f$ onto $F\cap y^\perp\,$, which is $f\,$. This implies
  $F\cap y^\perp\subset\varliminf_kE_k\,$, so $E_k\to F\cap y^\perp\,$. Moreover,
  \[
  x_k-E_k^*=x_k-F_k^*-\check E^*\to C-\check E^*=C\ ,
  \]
  so $\dom\vphi\subset F\cap y^\perp$ by the first part. Hence, $\dim\dom\vphi\sle\liminf_k\dim E_k\,$.

  We need to see that eventually, $y\not\perp F_k\,$. If it were true that $y\perp F_k$ frequently, then, passing to a subsequence,
  we could assume that $y\perp F_k$ \fa $k\,$. Hence, $y\perp F\,$. But this would imply $y\in\Span0F\cap F^\perp=0\,$,
  a contradiction.  Thus, eventually, $y\not\perp F_k\,$, so $\dim E_k\sle m-1\,$, and this proves that $\dim\dom\vphi<m\,$.
\end{Partsproof}

\begin{lem}\label{projection-closedgraph}
  The map
  \[
  (\pi,\lambda):Y_j\to P\times X^+
  \]
  has closed graph, where $X^+$ is the one point compactification of $X\,$.
\end{lem}

\begin{pf}
  Let $x_k-F_k^*\to x-F^*$ where $F_k,F\in P_j\,$ and $x_k\in F^\circledast\,,\,x\in F^\circledast\,$. Further, let $F_k\to E\in P\,$.
  By Lemma \ref{projection-subgraph-closed} (i), $F=\dom\sigma_{x-F^*}\subset E\,$. But $\dim E\sle n_{d-j}$ by
  Proposition \ref{lower-dim-closed}, which proves that $E=F\,$. By Lemma \ref{projection-subgraph-closed} (iii), $(x_k)$ is bounded, so by
  part (ii) of that lemma, $x=\limk_kx_k\,$.
\end{pf}

\par\noindent
Now we are ready to prove Theorem \ref{compact-bundle2}.

\begin{Partsproof}
\item[\emph{Proof of Theorem \ref{compact-bundle2} (i).}] If $\clos\Omega=U_{d+1}\,$, then the latter is compact. Let $F_k\in P\,$, \scth $F_k\to F\in\closed0X\,$. Then 
  $F$ is a convex cone. Passing to subsequences, we may assume
  $-F_k^*\to y-E^*\in U_{d+1}\,$. On the other hand, continuity of polarity gives $-F_k^*\to-F^*\,$. Thus
  $E=\dom\sigma_{y-E^*}=\dom\sigma_{-F^*}=F\in P\,$.\incrementcounter
\item The map $(\pi,\lambda):Y_j\to P\times X^+$ has compact range and closed graph, by Lemma \ref{projection-closedgraph}. Hence, it is continuous.
\item Let $x_k-F_k^*\to A\in\closed0X$ where $F_k\in P$ and $x_k\in F_k^\circledast\,$. Since $n=\dim X$ is finite,
  there exists $0\sle j\sle d$ such that $F_k\in P_j$ frequently. Passing to subsequences, we may assume $F_k\in P_j$ \fa $k\,$, and
  $F_k\to F\in P_j\,$. Let $C=\dom\sigma_A\,$. If $\dim C<n_{d-j}\,$, then $C\neq F\,$, and by
  Lemma \ref{projection-subgraph-closed} (ii), $(x_k)$ is unbounded. Lemma \ref{projection-subgraph-closed} (iii) provides us
  with $E_k,E\in P\,$, such that $x_k-E_k^*\to A\,$, $E_k\to E\,$, and $\dim E_k,\dim E<n_{d-j}\,$.

  We may write $x_k=u_k+v_k$ where $u_k\in E_k^\circledast$ and $v_k\perp E_k\,$. We claim that $u_k-E_k^*\to A\,$. Let $w\in A\,$. Then there exist $e_k\in E_k^*$
  such that $x_k-e_k\to w\,$. Then $v_k\in E_k^\perp\subset-E_k^*\,$, and
  \[
  w=\limk_k\Parens1{u_k+v_k-e_k}\in\varliminf\nolimits_k\Parens1{u_k-E_k^*}\ .
  \]
  Conversely, let $\alpha(\nats)$ be cofinal in $\nats$ and $e_{\alpha(k)}\in E_{\alpha(k)}^\circledast$ such that
  $w=\limk_k\Parens1{u_{\alpha(k)}-e_{\alpha(k)}}$ exists in $X\,$. Then $v_{\alpha(k)}\in E_{\alpha(k)}^\perp\subset E_k^*\,$, and
  \[
  w=\limk_k\Parens1{x_{\alpha(k)}-v_{\alpha(k)}-e_{\alpha(k)}}\in\limk_k\Parens1{x_k-E_k^*}=A\ .
  \]
  Thus, $\varlimsup_k\Parens1{u_k-E_k^*}\subset A\,$, and this proves our claim.

  Proceeding inductively (replace $x_k$ by $u_k$ and $F_k$ by $E_k$), we may assume that we have $\dim C=n_{d-j}\,$, \sth $(x_k)$ is bounded by
  Lemma \ref{projection-subgraph-closed} (iii). Then $C=F$ and $A=x-F^*$ where $x=\limk_kx_k\,$, by part (ii) of the lemma. Thus, we conclude 
  $\clos\Omega=U_{d+1}\,$, which proves the theorem.
\end{Partsproof}

\begin{cor}
  If the $P_j$ are compact \fa $j\,$, then
  \[
  \mathcal W_\Omega=\Set1{(x-F^*,y_1+y_2-x)}{x,y_1\in F^\circledast\,,\,y_2\in F^\perp\,,\,F\in P}\ ,
  \]
  with range and source are given by
  \[
  r(x-F^*,y_1+y_2-x)=x-F^*\nd s(x-F^*,y_1+y_2-x)=y_1-F^*
  \]
\end{cor}

\begin{pf}
  The condition $s(x-F^*,y)\in\clos\Omega$ reads $x+y\equiv y_1\pmod{F^\perp}$ \fs $y_1\in F^\circledast\,$, so we may set $y_2=x+y-y_1\,$.
\end{pf}

\subsection{Transversals and the Spectrum of $\mathrm C^*\bigl(\mathcal W_\Omega\bigr)$}
As is suggested by our study of $\clos\Omega\,$, we shall now \emph{always assume that $P_j$ be compact \fa $j\,$}. Moreover, somewhat abusing notation, we shall
identify $\clos\Omega$ with its image under $(\pi,\lambda)\,$, i.e.~we let $x-F^*\equiv(F,x)\,$. Of course, one should beware that the components of $(F,x)\in\clos\Omega$ 
depend  continuously on $x-F^*$ only when the latter is restricted to $Y_j\,$, and not globally.

  From \cite[ex.~2.7]{muhly-renault-williams-equivalences} recall that an \emph{abstract transversal} $T$ of some locally compact groupoid $\mathcal G$ is a closed subset
  of the unit space $\mathcal G^{(0)}$ meeting each orbit of the right action of $\mathcal G$ on $\mathcal G^{(0)}\,$, such that $r|\mathcal G_T$ and $s|\mathcal G_T$
  are open, where $\mathcal G_T=s^{-1}(T)\,$.

\begin{prop}\label{grp-equiv}
  The natural embedding $P\to\clos\Omega:F\mapsto(F,0)$ is a homeomorphism onto its closed image. Thus identified with its image in $\clos\Omega\,$, $P_j$ is an abstract
  transversal for $\mathcal W_\Omega|Y_j\,$.  Therefore, $M_j=s^{-1}(P_j)=\mathcal W_{\Omega,P_j}$ is a topological $(\mathcal W_\Omega|Y_j,\Sigma_j)$-equivalence,
  where $\Sigma_j=\mathcal W_\Omega|P_j$ is an Abelian group bundle with unit space $P_j\,$.
\end{prop}

\begin{pf}
  The continuity of the embedding is just continuity of polarity, cf.~\cite[cor.~7.2.12]{beer-topologies-book}. Since $P$ is compact, the embedding is topological. If $(F,x)\in Y_j\,$, 
  then the groupoid element $\gamma=(F,x,-x)\in\mathcal W_\Omega|Y_j$ satisfies $(F,x).\gamma=(F,0)\in P_j\,$, so $P_j$ meets every orbit in $Y_j\,$. To check the 
  openness of $r$ and $s$ on $M_j\,$, we first determine $M_j\,$. Indeed,
  \[
  M_j=\Set1{(F,x,y-x)}{F\in P_j\,,\,x\in F^\circledast\,,\,y\in F^\perp}\ .
  \]
  As to the openness of $r|M_j\,$, we can produce a section by
  \[
  \sigma:Y_j\to M_j:(F,x)\mapsto(F,x,-x)=(F,x,-\lambda(F,x))
  \]
  This section is continuous, since $\lambda$ is continuous on $Y_j\,$. Similarly, a section for $s|M_j$ is given by $\tau:P_j\to M_j:F\mapsto (F,0,0)\,$. This section actually
  extends continuously to a section of $s|\mathcal W_{\Omega,P}$ (which is not the case for $\sigma$). Thus $P_j$ is indeed an abstract transversal,
  and by \cite[ex.~2.7]{muhly-renault-williams-equivalences}, $M_j$ is therefore an equivalence.

  As to the last statement, it suffices to check that $r|\Sigma=s|\Sigma$ are trivial and that the groups $r^{-1}(F)=s^{-1}(F)$ are Abelian. To that end, note
  \[
  \Sigma_j=\mathcal W_\Omega|P_j=\Set1{(F,0,y)}{F\in P_j\,,\,y\perp F}\ ,
  \]
  so that $r$ and $s$ coincide, and their fibre at $F$ is $F^\perp\,$, with the usual group structure induced from the ambient vector space $X\,$. In passing, note that
  $\Sigma_j$ carries the relative topology induced from $P_j\times X\,$.
\end{pf}

\par\noindent
The existence of a Haar system for $\Sigma_j$ follows from the openness of its range projection, but can also be checked by hand as follows.

\begin{lem}
  The Abelian group bundle $\Sigma_j$ has a Haar system given by $\lambda^F=\delta_F\otimes\lambda_{F^\perp}\,$, where $\lambda_{F^\perp}$ denotes Lebesgue measure
  on the subspace $F^\perp\subset X$ endowed with the induced Euclidean structure.
\end{lem}

\begin{pf}
  We need to check the continuity. To that end, note that $F\mapsto F^\perp:P\to\closed0X$ is continuous by Proposition \ref{lower-dim-closed}. Let $m=n-n_{d-j}=\dim F^\perp$
  for $F\in P_j\,$. Fix $G\in P_j\,$. For some neighbourhood $U\subset P_j$ of $G\,$, $p_{F^\perp}:G^\perp\to F^\perp$ is an isomorphism \fa $F\in U\,$. For any
  $\vphi\in\Cc0X\,$,
  \[
  \int\vphi\,d\lambda_{F^\perp}=\int_{F^\perp}\vphi\,d\mathcal H^m
  =\sqrt{\det\Parens1{p_{F^\perp}^*p_{F^\perp}^{\vphantom*}}}\cdot\int\vphi\circ p_{F^\perp}^{-1}\,d\lambda_{G^\perp}\ ,
  \]
  by the area formula, cf.~\cite[cor.~3.2.20]{federer-geommeas}. (Here, $\mathcal H^m$ denotes $m$-dimensional Hausdorff measure.) The continuity follows from Lebesgue's 
  dominated convergence theorem.
\end{pf}

\begin{cor}\label{morita-equiv}
  There is a completion of $\Cc0{M_j}$ to a $\Parens1{\Cst0{\mathcal W_\Omega|Y_j},\Cst0{\Sigma_j}}$ equivalence bimodule $\Cst0{M_j}\,$, thus establishing a strong
  Morita equivalence $\Cst0{\mathcal W_\Omega|Y_j}\sim\Cst0{\Sigma_j}\,$. Moreover, $\Cst0{\Sigma_j}\cong\Ct[_0]0{\Sigma_j}$ by Fourier transform.
  In particular, $\Cst0{\mathcal W_\Omega|Y_j}$ is liminary, of spectrum $\Sigma_j\,$.
\end{cor}

\begin{pf}
  Strong Morita equivalence follows from Proposition \ref{grp-equiv} and \cite[th.~2.8]{muhly-renault-williams-equivalences}. Define, for $\vphi\in\Cc0{\Sigma_j}\,$, the
  fibrewise Fourier transform
  \[
  \mathcal F(\vphi)(F,y)=\int_{F^\perp}e^{-2\pi i\Rscp0y\eta}\vphi(F,\eta)\,d\eta\ .
  \]
  From Euclidean Fourier analysis, $\mathcal F:(\Cc0{\Sigma_j},*)\to\Ct[_0]0{\Sigma_j}\,$, is a continuous $*$-morphism. Thus, there exists an extension to a
  contractive $*$-mor\-phism $\Cst0{\Sigma_j}\to\Ct[_0]0{\Sigma_j}\,$. The image of $\mathcal F$ is dense by Stone--Weierstrass. Indeed, two points
  $(F_1,y_1),(F_2,y_2)\in\Sigma_j$ with $F_1\neq F_2$ are easily separated. If $F_1=F_2=F$ and $y_1\neq y_2\,$, we can separate $y_1$ and $y_2$ by the
  Fourier transform on $F^\perp$ of some $\vphi\in\Cc0{F^\perp}\,$. Now consider
  \[
  \psi(F',y)=\chi(F')\cdot\vphi\Parens1{p_{F^\perp}(y)}\mathfa (F',y)\in\Sigma_j\ ,
  \]
  where $\chi\in\Ct0{P_j}$ \scth $\chi(F)=1\,$. $\psi$ is continuous because the $P_j$ are compact, and by Corollary \ref{proj-cont} (ii). Moreover, $\mathcal F\psi$ separates
  $(F,y_j)\,$, $j=1,2\,$.

  Thus, clearly, the locally compact space $\Sigma_j$ injects into the spectrum of the commutative C$^*$-algebra $\Cst0{\Sigma_j}\,$. Conversely, let $\chi$ be a
  character of $\Cst0{\Sigma_j}\,$. Denote by $I_F$ the ideal of $\Cst0{\Sigma_j}$ generated by the functions vanishing on the fibre of $\Sigma_j$ over $F\,$.
  A partition of unity argument shows that for $F\neq F'\,$, $F,F'\in P_j\,$, $I_F+I_{F'}=\Cst0{\Sigma_j}\,$. Thus, there exists a unique $F\in P$ for which
  $\chi(I_F)\neq0\,$. Clearly, $\Cst0{F^\perp}=\Cst0{\Sigma_j}/I_F\,$, so $\chi$ is given by the Fourier transform with respect to $F^\perp\,$,
  evaluated at some $y\in F^\perp\,$.

  Hence, $\Sigma_j$ exhausts the spectrum, and Gelfand's theorem shows that $\mathcal F$ is injective, and therefore an isometric $*$-isomorphism.

  Since strong Morita equivalence implies stable isomorphism for $\sigma$-unital C$^*$-algebras, we find that $\Cst0{\mathcal W_\Omega}$ is stably isomorphic to
  $\Ct[_0]0{\Sigma_j}\,$. Indeed, the separability of these C$^*$-algebras follows from the separability of their underlying groupoids.
\end{pf}

\begin{rem}\ 
\begin{enumerate}\romannum
\item We have $\Sigma_d=\{0\}\times X=Y_d\,$, so $\Cst0{\mathcal W_\Omega|Y_d}\cong\Ct[_0]0X\,$. Similarly, $\Sigma_0=\{\Omega^*\}\times0$ and $Y_0=\Omega\,$, so
  $\Cst0{\mathcal W_\Omega|\Omega}\cong\mathbb K\,$.
\item As follows by the theory of $*$-algebraic bundles, the latter statement of the above corollary is true for any Abelian group bundle endowed with a Haar system,
  cf.~\cite[th.~1.3.3]{ramazan-thesis}.
\item There is a delicate point to the above equivalences. Namely, although $P$ itself is a compact subset of the unit space $\clos\Omega\,$, meeting
each orbit, it is not an abstract transversal in the sense defined above. Indeed, $\Cst0{\mathcal W_\Omega}$ has a faithful irreducible representation, so its spectrum
contains a dense point, and unless $X=0\,$, the spectrum is non-Hausdorff. However, if $P$ were a transversal, then
$\Cst0{\mathcal W_\Omega}$ would be Morita equivalent to a commutative C$^*$-algebra.

\begin{center}
\includegraphics[scale=0.3, clip=true]{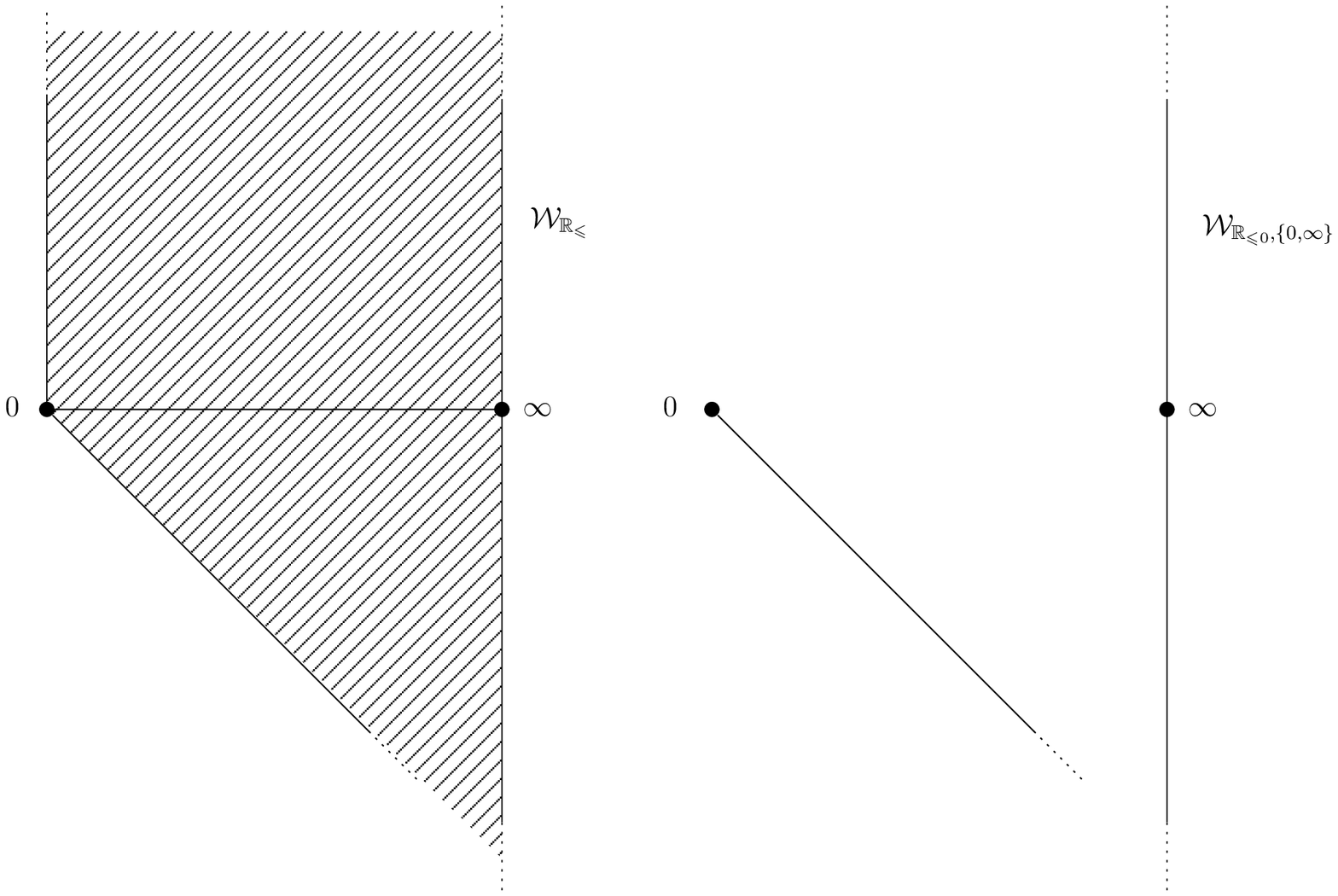}
\end{center}

Which condition fails can be inspected for the example of the classical Wiener--Hopf algebra where $X=\reals$ and $\Omega=\reals_{\sge0}\,$. 
Then $\clos X=[-\infty,\infty]$ under the identification $x\mapsto\Lopint0{-\infty}x\,$, and this gives the order topology for this interval. Similarly,
$\clos\Omega=[0,\infty]\,$. In this realisation, the action of $\reals$ is by translation on $\reals$ and trivial at $\pm\infty\,$. Thus, $\mathcal W_{\reals_{\sge0}}$ and
$\mathcal W_{\reals_{\sge0},\{0,\infty\}}$ (note $P=\{0,\infty\}$) work out as in the illustration. The range and source projection are given by 
\[
r(x,-x)=x\,,\,r(\infty,x)=\infty\nd s(x,-x)=0\,,\,s(\infty,x)=\infty\ . 
\]
Although $s$ is open (it always is), $r$ is not, since an open neighbourhood of $\infty$ projects to the non-open point $\infty\in[0,\infty]\,$. The first named author wishes 
to thank George Skandalis for pointing out this observation.
\end{enumerate}
\end{rem}

\begin{thm}
  The sets $U_j\subset\clos\Omega\,$, $j=0,\dotsc,d+1\,$, form an ascending chain of open invariant subsets. The ideals $I_j=\Cst0{\mathcal W_\Omega|U_j}$ form a
  composition series with liminary quotients $I_{j+1}/I_j=\Cst0{\mathcal W_\Omega|Y_j}\,$. Hence, the C$^*$-algebra $\Cst0{\mathcal W_\Omega}$ is of type I
  (i.e.~postliminary). Its spectrum is the set $\Sigma=\bigcup_{j=0}^d\Sigma_j\,$, the topology of which is given by the sets $U\cup\bigcup_{i=0}^{j-1}\Sigma_i$ \fa
  $0\sle j\sle d$ and all open $U\subset\Sigma_j\,$.
\end{thm}

\begin{pf}
  The computation of the quotients is given by any of the following sources: \cite[ch.~II, prop.~4.5]{renault-groupoid-cst}, \cite[2.4]{hilsum-skandalis},
  \cite[prop.~2.4.2]{ramazan-thesis}. We already know $\Cst0{\mathcal W_\Omega|Y_j}\sim\Ct[_0]0{\Sigma_j}\,$, so this algebra is liminary of spectrum $\Sigma_j\,$.

  Set $\Sigma=\widehat{\Cst0{\mathcal W_\Omega}}\,$. Let
  \[
  V_j=\Set1{\vrho\in\Sigma}{\vrho(I_{j+1})\neq0}\nd W_j=\Set1{\vrho\in\Sigma}{\vrho(I_j)=0}\ .
  \]
  Then $V_j$ is open, and $W_j$ is closed, and we have $V_j\cap W_j\approx\widehat{I_{j+1}/I_j}\approx\Sigma_j$ \cite[prop.~3.2.1]{dixmier}. If $U\subset\Sigma_j$ is open,
  then $\Sigma_j\setminus U$ is closed in $V_j\,$, and hence $V_j\setminus(\Sigma_j\setminus U)=U\cup\bigcup_{i=0}^{j-1}\Sigma_i$ is open in $\Sigma\,$.

  Conversely, since $\Cst0{\mathcal W_\Omega}$ has a faithful unitary representation, $V_0=V_0\cap W_0\approx\Sigma_0=*$ is dense in $\Sigma\,$. Hence, any open
  $\vvoid\neq V\subset\Sigma$ is dense. The assertion follows.
%
%
\end{pf}

\begin{cor}
  The C$^*$-algebra $\Cst0{\mathcal W_\Omega}$ is solvable of length $d\,$, in the sense of \cite{dynin-invsing}.
\end{cor}

\begin{pf}
  The above composition series is uniquely determined by the requirement that $I_{j+1}/I_j$ be the largest liminary subalgebra of $I_{d+1}/I_j\,$, by \cite[prop.~4.3.3]{dixmier}.
  That this requirement obtains in turn follows from \cite[prop.~4.2.6]{dixmier}. Hence, the length is exactly $d\,$.
\end{pf}

\section{Analytical Indices}

\subsection{Continuous Fields of Elementary C$^*$-Algebras}\label{jth-wh-ext}
  The above considerations show that we have short exact sequences
  \[
  \begin{CD}
    0@>>>{\Cst0{\mathcal W_\Omega|Y_{j-1}}}@>>>I_{j+1}/I_{j-1}@>>>{\Cst0{\mathcal W_\Omega|Y_j}}@>>>0\ .\end{CD}
  \]
  These may be considered as elements $\partial_j\in KK^1\Parens1{\Ct[_0]0{\Sigma_j},\Ct[_0]0{\Sigma_{j-1}}}\,$, and the corresponding homomorphisms of the $K$-groups 
  are then given by the Kasparov product with $\partial_j\,$. In order to give an analytical description of the $\partial_j\,$, we have to compute the subquotients 
  $\Cst0{\mathcal W_\Omega|Y_j}$ of the composition series more explicitly. In fact, we shall exhibit them as continuous fields of elementary C$^*$-algebras, thereby giving 
  an independent proof of results from section $2\,$.\vspace{1ex}

  \par\noindent
  Fix $0\sle j\sle d$ and $E\in P_j\,$. A pair $(\psi_U,U)$ where $U\subset P_j$ is an open neighbourhood of $E$ and $\psi_U:U\times X\to X$ shall be called a
  \emph{positive local trivialisation at $E$} if the following conditions are satisfied:
  \begin{enumerate}\romannum\parindent0pt
  \item $\psi_U$ is continuous, $\psi_F=\psi_U(F,\var):X\to X$ is bi-Lipschitz \fa $F\in U\,$,
  \item for a.e.~$x\in F^*\,$, $\psi_F'(x)$ exists, and $\det\psi_F'(x)>0\,$,
  \item \fa $F\in U\,$, $\psi_F$ is linear when restricted to $F^\perp\,$, and
  \item \fa $F\in U\,$, $\psi_F(F^\circledast)=E^\circledast$ and $\psi_F(F^\perp)=E^\perp\,$.
  \end{enumerate}
  We point out that the derivative exists for almost every $x\in X\,$, by Rademacher's theorem, cf.~\cite[th.~3.1.6]{federer-geommeas}. If, moreover,
  $\det\psi_F'(x)=1$ whenever $\psi_F'(x)$ exists, then the local trivialisation shall be termed \emph{normalised}.

\begin{prop}\label{groupoid-loc-triv}
  Let $(\psi_U,U)$ be a normalised local trivialisation $E\in P_j\,$. Then there exists a $*$-isomorphism
  \[
  \Psi_U:\Cst0{\mathcal W_\Omega|\pi^{-1}(U)}\to\Ct[_0]0{U\times E^\perp}\otimes\Cst0{\mathcal W_{E^\circledast}|E^\circledast}
  \]
  given by
  \[
  \Psi_U(\vphi)(F,y,u,v-u)=\int_{E^\perp}e^{-2\pi i\Rscp0yz}\vphi\Parens1{F,\psi_F^{-1}(u),\psi_F^{-1}(v)+\psi_F^{-1}(z)-\psi_F^{-1}(u)}\,dz
  \]
  \fa $\vphi\in\Cc0{\mathcal W_\Omega|\pi^{-1}(U)}\,$, $F\in U\,$, $y\in E^\perp\,$, $u,v\in E^\circledast\,$. Here, recall that 
  \[
  \mathcal W_{E^\circledast}|E^\circledast=(\Span0E\rtimes\Span0E)|E^\circledast\ .
  \]
\end{prop}

\begin{pf}
  We have seen in Corollary \ref{morita-equiv} how fibrewise Fourier transform establishes an isomorphism $\Cst0{U\times E^\perp}\cong\Ct[_0]0{U\times E^\perp}\,$, where
  $U\times E^\perp$ is considered an Abelian group bundle.

  Thus, it remains to see that
  \[
  \Phi_U:\left\{
    \begin{array}[c]{l}
      \mathcal W_\Omega|\pi^{-1}(U)\to \Parens1{U\times E^\perp}\times\Parens1{\mathcal W_{E^\circledast}|E^\circledast}\\
      \Parens1{F,x,y_1+y_2-x}\mapsto\Parens1{F,\psi_F(y_2),\psi_F(x),\psi_F(y_1)-\psi_F(x)}
    \end{array}\right.
  \]
  is a topological isomorphism of groupoids in the sense of \cite{muhly-renault}. That it is a groupoid isomorphism is clear from the linearity of $\psi_F$ on $F^\perp\,$; 
  moreover, it is immediate that it is a homeomorphism. Finally, the Haar system of $\mathcal W_\Omega$ is $\lambda^{F,x}=\delta_{(F,x)}\otimes\lambda_{F^*-x}\,$.
  We find, by the change of variables formula,
  \[
  \Phi_U(\lambda^{F,x})=\Phi_U(\delta_{(F,x)}\otimes\lambda_{F^\perp}\otimes\lambda_{F^\circledast-x})
  =\Abs0{\det\psi_F'}\cdot\Parens1{\delta_F\otimes\lambda_E^\perp\otimes\delta_{\psi_F(x)}\otimes\lambda_{E^\circledast-\psi_F(x)}}\ ,
  \]
  as required, since $\det\psi_F'=1$ a.e.
\end{pf}

\par\noindent
  In order to make this proposition substantial, we need to construct normalised local trivialisations. It is clear that given a positive local trivialisation, it can be normalised.
  Moreover, Corollary \ref{proj-cont} shows that for $F$ close to $E\,$, $p_{F^\perp}$ is a linear isomorphism of $F^\perp$ onto $E^\perp\,$. Hence, $\psi=\psi_U$ can be
  constructed as $\psi(F,x)=p_E\psi_1(F,p_F(x))+p_{E^\perp}y$ as soon as a map $\psi_1$ can be given which satisfies all the conditions of a local trivialisation, apart from
  linearity on $F^\perp$ and $\psi_F(F^\perp)=E^\perp\,$.

\begin{prop}\label{norm-loc-triv-ex}
  For $E\in P_j\,$, there exists an open neighbourhood $E\in U\subset P_j$ and a map $\psi_U:U\times X\to X$ such that
  \begin{enumerate}\romannum\parindent0pt
  \item $\psi_U$ is continuous, \fa $F\in U\,$, $\psi_F=\psi_U(F,\var)$ is bi-Lipschitz,
  \item for a.e.~$x\in F^*\,$, the derivative $\psi_F'(x)$ exists and $\det\psi_F'(x)>0\,$, and
  \item $\psi_F(F^\circledast)=E^\circledast\,$.
  \end{enumerate}
  In particular, there exists a normalised local trivialisation at $E\,$.
\end{prop}

  First, note the following definition and lemma. Any $x\in\partial C$ where $C\subset X$ is closed and convex
  with $C^\circ\neq\vvoid\,$, is called a $C^1$-point, if there is a unique supporting hyperplane at $x\,$.

\begin{lem}\label{minkowski-directional-deriv}
  Let $C\subset X$ be a compact convex neighbourhood of $0\,$, and $\mu:X\to\Ropint00\infty$ denote its Minkowski gauge functional, i.e.
  \[
  \mu(x)=\inf\Set1{\alpha>0}{\alpha^{-1}\cdot x\in C}\ .
  \]
  \Fa $x\in X\setminus0\,$, $v\in X\,$, the right directional derivative $\nabla_v^+\mu(x)=\tfrac d{dt}\mu(x+tv)\bigm\vert_{t=0+}$ exists, and
  \[
  \nabla_v^+\mu(x)=\sigma_{n_x(C)}(v)\mathtxt{where}n_x(C)=\Set1{y\in N_{\mu(x)^{-1}\cdot x}(C)}{\Rscp0xy=\mu(x)}
  \]
  and $N_x(C)=\Set1{y\in X}{\Rscp0yx=\sigma_C(y)}$ is the normal cone of $C$ at $x\,$. In particular, $\mu$ is differentiable at $x\in X\setminus0$ if and only if
  $\mu(x)^{-1}\cdot x$ is $C^1\,$, with gradient
  \[
  \nabla\mu(x)=\frac{\mu(x)}{\Rscp1{\pi_{N_{\mu(x)^{-1}x}(C)}(x)}x}\cdot\pi_{N_{\mu(x)^{-1}\cdot x}(C)}(x)\ .
  \]
  Here, we remind the reader that $\pi_C$ denotes the metric projection onto the closed convex set $C\,$, cf.~\cite{zarantonello-convex-spectral}. 
\end{lem}

\begin{pf}
  Since $\mu$ is convex, \cite[\S~3.2.(i), th.1]{giles-convan-diff} the right directional derivative $\nabla_v^+\mu(x)$ exists
  everywhere and defines a sublinear functional in $v\,$. By the Hahn--Banach theorem, $\nabla_v^+\mu(x)$ is the upper envelope of
  linear functionals it majorises. By \cite[\S~3.2.(i), th.3]{giles-convan-diff}, the subdifferential of $\mu$ at $x$ is
  \[
  \partial\mu(x)=\Set1{y\in X}{-\nabla_{-v}^+\mu(x)\sle\Rscp0yv\sle\nabla_v^+\mu(x)}\ .
  \]
  Since indeed $-\nabla_{-v}^+\mu(x)\sle\nabla_v^+\mu(x)\,$, we find
  \[
  \nabla_v^+\mu(x)=\sup\Set1{\Rscp0yv}{-\nabla_{-v}^+\mu(x)\sle\Rscp0yv\sle\nabla_v^+\mu(x)}
  =\sup\nolimits_{y\in\partial\mu(x)}\Rscp0yv\ .
  \]
  Since $\mu$ is positively $1$-homogeneous, $\nabla_v^+\mu(x)=\nabla_v^+\mu(rx)$ \fa $r>0\,$. Hence, we may restrict attention to
  the case $x\in\partial C\,$, i.e.~$\mu(x)=1\,$. By \cite[\S~3.2.(i), lem.~to th.~4]{giles-convan-diff}, we have
  \[
  \partial\mu(x)=\Set1{y\in X}{1=\Rscp0yx\sge\Rscp0yz\smathfa z\in C}=\Set1{y\in X}{1=\Rscp0yx=\sigma_C(y)}\ .
  \]
  Thus,
  \[
  \nabla_v^+\mu(x)=\sup\Set1{\Rscp0yv}{1=\Rscp0yx=\sigma_C(y)}=\sigma_{n_x(C)}(v)\ .
  \]
  This proves the first assertion.

  As to the second, if $y=\mu(x)^{-1}\cdot x$ is a $C^1$-point, the normal cone is just the ray $\reals_{\sge0}\cdot\pi_{N_y(C)}(x)\,$. Note
  \[
  \Rscp1x{r\cdot \pi_{N_y(C)}(x)}=\mu(x)\ \equiva\ r=\mu(x)\cdot\Rscp1{\pi_{N_y(C)}(x)}x^{-1}\ ,
  \]
  which implies
  \[
  \nabla_v^+\mu(x)=\frac{\mu(x)}{\Rscp1{\pi_{N_y(C)}(x)}x}\cdot\Rscp1{\pi_{N_x(C)}(x)}v\ .
  \]
  Since the $C^1$-points of $C$ are exactly the boundary points of $C$ at which $\mu$ is differentiable, by
  \cite[\S~3.2.(i), th.~5]{giles-convan-diff}, the assertion follows.
\end{pf}

\begin{rem}
  We note that the above formula for $\nabla\mu(x)$ \emph{at $C^1$-points} also follows from \cite[prop.~3.1]{fitzpatrick-phelps}, who prove
  \[
  \nabla\mu(x)=\frac{x-\pi_C(x)}{\Rscp0{\pi_C(x)}{x-\pi_C(x)}}\mathfa x\in X\setminus C
  \]
  at which $\mu$ is differentiable. 

  Indeed, let $x\in\partial C\,$, $y=\pi_{N_x(C)}(x)$ and $z=\mu(y)^{-1}\cdot y\,$. Then $N_z(C)=N_x(C)\,$, and moreover, $\pi_C^{-1}(z)=z+N_z(C)\,$, by
  \cite[\S~2, lem.~2.4]{zarantonello-convex-spectral}. We have $\pi_C\Parens1{r\cdot y}=z$ \fa $r\sge\mu(y)^{-1}\,$, if $x$ is a $C^1$-point. This implies the above formula.
\end{rem}

{\def\Elproofname{\protect{PROOF of Proposition \ref{norm-loc-triv-ex}.}}
\begin{pf}
  We are done once we have constructed a map satisfying (i) and (ii), and which maps $p_{E^\perp}(F^\circledast)$ to $E^\circledast\,$, since $p_E$ sets up an isomorphism
  $\Span0F\to\Span0E$ for $F$ close to $E\,$, by Corollary \ref{proj-cont}. So we may as well assume that $F\subset\Span0E$ for all $F\in U\,$. Since we may then let $\psi_F$
  be the identity on $E^\perp\,$, for simplicity, we may assume the cones we are considering to be solid in $X\,$. 

  Let $\xi_0\in E^\circ\cap E^{*\circ}\,$, $\Norm0{\xi_0}=1\,$, and $H=\Set0x{\Rscp0x{\xi_0}=1}\,$. For $F$ close to $E\,$, we have $\xi_0\in F^\circ\cap F^{*\circ}\,$, too.
  Take $X_+=\xi_0^{*\circ}\subset X\setminus\xi_0^\perp$ to be the half-space containing $E^*\setminus0\,$. Then $F^*\subset X_+$ for $F$ close
  to $E\,$.  Let
  \[
  \mu_F(x)=\inf\Set1{\alpha>0}{\alpha^{-1}x\in H\cap F^*-\xi_0}\mathfa x\in\xi_0^\perp=H-\xi_0\ ,
  \]
  the Minkowski functional of the compact convex set $C_F=H\cap F^*-\xi_0\,$, which is a neighbourhood of zero in $\xi_0^\perp\,$. Let
  \[
  \vphi_F(x)=\frac{\mu_F(x)}{\mu_E(x)}\cdot x\mathfa x\in \xi_0^\perp\ .
  \]
  Then $\vphi_F:\xi_0^\perp\to\xi_0^\perp\,$, mapping $C_F$ to $C_E\,$. We may now define
  \[
  \psi_F\Parens1{x+r\cdot\xi_0}=\vphi(x)+r\cdot\xi_0\mathfa x\perp\xi_0\,,\,r\in\reals\ .
  \]
  In particular, $\psi_F=\vphi_F$ on $\xi_0^\perp\,$, and
  \[
  \psi_F(x)=\Rscp0x{\xi_0}\cdot\Parens2{\vphi_F\Parens2{\frac x{\Rscp0x{\xi_0}}-\xi_0}+\xi_0}\mathfa x\in X_+\ .
  \]
  Then condition (iii) is clearly verified.

  As to condition (i), we may assume $\ball(r,\xi_0)\subset H\cap F^*\subset\ball(R,\xi_0)$ \fa $F$ and some $0<r<R\,$. This implies
  $r\cdot\Norm0\var\sle\mu_F\sle R\cdot\Norm0\var\,$. Assume that $\mu_E(x)\sge\mu_E(y)\,$. Then
  \[
  \mu_E\Parens1{\vphi(x)-\vphi(y)}\sle\mu_F(x)\cdot\mu_E\Parens2{\frac x{\mu_E(x)}-\frac y{\mu_E(y)}}
  +\mu_E\Parens2{\Parens1{\mu_F(x)-\mu_F(y)}\cdot\frac y{\mu_E(y)}}
  \]
  Further,
  \[
  \mu_E\Parens3{\Parens1{\mu_F(x)-\mu_F(y)}\cdot\frac y{\mu_E(y)}}\sle\frac1r\cdot \mu_F(x-y)\sle\frac Rr\cdot\Norm0{x-y}\ ,
  \]
  and
  \begin{align*}
    \mu_E\Parens3{\frac x{\mu_E(x)}-\frac y{\mu_E(y)}}&\sle\frac1{\mu_E(x)}\cdot\mu_E(x-y)+\mu_E\Parens1{\Parens1{\mu_E(x)^{-1}-\mu_E(y)^{-1}}\cdot y}\\
    &\sle\frac R{\mu_E(x)}\cdot\Norm0{x-y}+\frac Rr\cdot\Abs0{\mu_E(x)^{-1}-\mu_E(y)^{-1}}\cdot\mu_E(y)\\
    &=\frac R{\mu_E(x)}\cdot\Norm0{x-y}+\frac Rr\cdot\frac{\mu_E(x)-\mu_E(y)}{\mu_E(x)}\\
    &\sle\frac R{\mu_E(x)}\cdot\Parens2{1+\frac Rr}\cdot\Norm0{x-y}\ .
  \end{align*}
  Exchanging the role of $x$ and $y\,$, and noting that $\tfrac{\mu_F}{\mu_E}\sle\tfrac Rr\,$, we find that $\vphi$ is $L$-Lipschitz, where
  \[
  L=\frac R{r^2}\cdot\Parens2{1+R\Parens2{1+\frac Rr}}\ .
  \]
  It follows that $\psi_F$ is $L'$-Lipschitz with $L'=\sqrt 2\cdot\max(L,\Norm0{\xi_0})\,$. Since $\psi_F^{-1}$ is given by exchanging the roles of $E$ and $F$ in the
  definition of $\vphi_F\,$, it follows that $\psi_F$ is bi-Lipschitz. As to the joint continuity of $\psi\,$, it suffices to note
  \[
  \Norm1{\vphi_{F_1}(x)-\vphi_{F_2}(x)}=\Abs0{\mu_{F_1}(x)-\mu_{F_2}(x)}\cdot\Norm1{\mu_E(x)^{-1}\cdot x}\sle R\cdot\Abs0{\mu_{F_1}(x)-\mu_{F_2}(x)}\ ,
  \]
  and that $\mu_F$ depends continuously on $F\,$.

  Suffices to compute derivatives on $X_+\,$. By Lemma \ref{minkowski-directional-deriv} for $x\perp\xi_0\,$,
  \[
  \nabla\mu_F(x)=\frac{\mu_F(x)}{\Rscp0{\pi_{x,F}(x)}x}\cdot \pi_{x,F}(x)\mathtxt{whenever the derivative exists,}
  \]
  $\pi_{x,F}$ denoting the metric projection onto the normal cone $N_{\mu(x)^{-1}\cdot x}(C_F)\,$.

  A simple calculation gives \fa $x\perp\xi_0$ for which the derivative exists,
  \[
  \vphi_F'(x)v=\lambda\cdot v+\lambda\cdot\vrho_x(v)\cdot x\mathtxt{where}
  \vrho_x(v)=\frac{\Rscp0{\pi_{x,F}(x)}v}{\Rscp0{\pi_{x,F}(x)}x}-\frac{\Rscp0{\pi_{x,E}(x)}v}{\Rscp0{\pi_{x,E}(x)}x}
  \]
  and $\lambda=\frac{\mu_F(x)}{\mu_E(x)}>0\,$. Let $\xi_1=\Norm0x^{-1}\cdot x\,$,
  and complete this to an orthonormal basis $\xi_1,\dotsc,\xi_{n-1}$ of $\xi_0^\perp\,$. Then since $\vrho_x(x)=0\,$, $\vphi_F'(x)$ has the matrix expression
  \[
  \vphi_F'(x)=\begin{Matrix}1\lambda& 0 & \dotsm & 0\\ \vrho_x(\xi_2) & \lambda & & \\ \vdots & & \ddots & \\ \vrho_x(\xi_{n-1}) & & & \lambda\end{Matrix}\ .
  \]
  We find $\det\vphi_F'(x)=\lambda^{n-1}>0\,$.

  Next, let $x\in X_+=\reals_{>0}\cdot\xi_0\oplus\xi_0^\perp$ be arbitrary, and define $N(x)=\tfrac x{\Rscp0x{\xi_0}}\,$. Then
  \[
  \psi_F'(x)v=\Rscp0v{\xi_0}\cdot\Parens1{\vphi_F(N(x)-\xi_0)+\xi_0}+\Rscp0x{\xi_0}\cdot\vphi_F'(N(x)-\xi_0)N'(x)v\ .
  \]
  Let $\lambda=\tfrac{\mu_F}{\mu_E}\Parens1{N(x)-\xi_0}$ and $\xi_1=\Norm0{N(x)-\xi_0}^{-1}\cdot\Parens1{N(x)-\xi_0}\,$. Observe
  \[
  \Rscp0x{\xi_0}N'(x)v=v-\frac{\Rscp0v{\xi_0}}{\Rscp0x{\xi_0}}\cdot x=
  \begin{cases}
    \Rscp0v{\xi_0}\cdot\Parens1{\xi_0-N(x)} & v\in\reals\cdot\xi_0\ ,\\
    v & v\perp\xi_0\ .
  \end{cases}\ .
  \]
  In particular, we note that $N'(x)v\perp\xi_0$ for every $x\,$, and since $\vphi_F$ is $1$-homogeneous,
  \[
  \vphi_F'\Parens1{N(x)-\xi_0}(\xi_0-N(x))=-\vphi_F(N(x)-\xi_0)
  \]
  In terms of the orthonormal basis $\xi_0,\xi_1,\dotsc,\xi_{n-1}\,$, $\psi_F'(x)$ has the matrix expression
  \[
  \psi_F'(x)=\begin{Matrix}1
    1 & 0 & \dotsm & 0\\
    0 & & & \\
    \vdots & & \vphi_F'(N(x)-\xi_0)  & \\
    0 & & & \\
  \end{Matrix}\ .
  \]
  In particular,
  \[
  \det\psi_F'(x)=\det\vphi_F'(x)=\lambda^{n-1}=\frac{\mu_F}{\mu_E}\Parens1{N(x)-x_0}^{n-1}>0\ .
  \]
  This proves the proposition.
\end{pf}}

\begin{cor}
  For $0\sle j\sle d\,$, $M_j$ is a oriented real vector bundle over $Y_j$ of rank $n-n_{d-j}\,$. Similarly, $\Sigma_j$ is an oriented real vector bundle over $P_j$ of
  rank $n-n_{d-j}\,$.
\end{cor}

\par\noindent
  Let $E\in P_j\,$. Then, for the groupoid $\mathcal G=\mathcal W_{E^\circledast}|E^\circledast\,$, we have
  \[
  \mathcal G_v=s^{-1}(v)=\Set1{(u,v-u)}{u\in E^\circledast}\mathfa v\in E^\circledast\ ,
  \]
  so we may identify $\Lp[^2]0{\mathcal G_v}$ with $\Lp[^2]0{E^\circledast}\,$. The regular representation $\vrho_E$ of $\Cst0{\mathcal W_{E^\circledast}|E^\circledast}$
  on this space is given by
  \[
  \vphi*h(u)=\int_{E^\circledast}\vphi(u,w-u)h(w)\,dw
  \mathfa\vphi\in\Cc0{\mathcal W_{E^\circledast}|E^\circledast}\,,\,h\in\Lp[^2]0{E^\circledast}\ .
  \]
  This is manifestly independent of $v\,$. In the notation of \cite[2.12.1-2]{muhly-renault}, the representation $\vrho_E$ is just $J^{-1}\ind\delta_0\,J\,$.

  On the other hand, for $(F,y)\in\Sigma_j$ define $L_\Omega^{F,y}=L_{\vphantom\Omega}^{F,y}$ by
  \begin{align*}
    L^{F,y}(\vphi)h(v)&=\int_{F^*-v}\vphi(F,v,w)e^{-2\pi i\Rscp0wy}h\Parens1{v+p_F(w)}\,dw\\
    &=\int_{F^\perp}\!\int_{F^\circledast}\vphi(F,v,w_1+w_2-v)e^{-2\pi i\Rscp0{w_2}y}h(w_1)\,dw_1\,dw_2
  \end{align*}
  \fa $\vphi\in\Cc0{\mathcal W_\Omega}\,$, $h\in\Lp[^2]0{F^\circledast}\,$, and $v\in F^\circledast\,$. The following proposition is then straightforward.

\begin{prop}\label{equiv-rep}
  Let $E\in P_j\,$, $(\psi_U,U)$ a normalised local trivialisation, and fix some $(F,y)\in\Sigma_j|U\,$. If $\chi_{F,y}$ denotes the character of
  $\Ct[_0]0{U\times E^\perp}$ given by evaluation at $(F,z)\,$, where $(\psi_F|F^\perp)^tz=y\,$, then $\Parens1{\chi_{F,y}\otimes\vrho_E}\circ\Psi_U$ and
  $L^{F,y}$ are equivalent representations of $\Cst0{\mathcal W_\Omega|\pi^{-1}(U)}\,$.
\end{prop}

\par\noindent
  For measurable $E\subset X$ and functions $f,g\,$, define the following abbreviations whenever they make sense:
  \begin{gather*}
    \mathcal F_E(f)(x)=\int_Ee^{-2\pi i\Rscp0xy}f(y)\,d\lambda_{\Span0E}(y)\ ,\ \mathcal F_E^*(f)=\int_Ee^{2\pi i\Rscp0xy}f(y)\,d\lambda_{\Span0E}(y)\ ,\\
    f*_Eg(x)=\int_Ef(y)g(x-y)\,d\lambda_{\Span0E}(y)\ ,\\
    f^*(x)=\overline{f(-x)}\ ,\ f^y(x)=f(x+y)\ ,\ e_y(x)=e^{2\pi i\Rscp0xy}\ .
  \end{gather*}
  Note the following equations:
  \begin{gather*}
    (f*_Eg)*_Fh=f*_E(g*_Fh)\smathtxt{for}F-E=F\ ,\ f*_Eg(x)=g*_{x-E}f(x)\smathtxt{for}x\in E\ ,\\
    f^**_Eg^*(x)=(g*_{E-x}f)^*(x)\smathtxt{for}x\in\Span0E\ ,\ f^y*_Eg=f*_{E+y}g^y\smathtxt{for}y\in\Span0E\ ,\\
    \mathcal F_F(f)*_E\mathcal F_G(g)=\mathcal F_{F\cap G}(f\cdot g)\smathtxt{for}E=\Span0F=\Span0G=\Span0{F\cap G}\ ,\\
    \mathcal F_E(f^y)=e_y\cdot\mathcal F_{E+y}(f)\smathtxt{for}y\in\Span0E\ ,
  \end{gather*}
  which are standard applications of Euclidean Fourier analysis.

\begin{prop}\label{cont-field}
  The family $\mathcal E_j=(\Lp[^2]0{F^\circledast})_{(F,y)\in\Sigma_j}$ is a continuous field of Hilbert spaces with a dense subspace $\Theta$ of sections given by the maps
  $(F,y)\mapsto\vphi_{F,y}\,$, where
  \[
  \vphi_{F,y}(x)=\mathcal F_{F^\perp}(\vphi^x)(y)\mathfa\vphi\in\Cc0X\,,\,(F,y)\in\Sigma_j\,,\,x\in F^\circledast\ .
  \]
\end{prop}

\begin{pf}
  By \cite[prop.~10.2.3]{dixmier}, it suffices to show that $\Theta$ is dense in every fibre and that $\Norm0\vphi$ is continuous \fa $\vphi\in\Theta\,$. The density follows 
  by considering the algebraic tensor product $\Cc0{F^\circledast}\odot\Cc0{F^\perp}\,$. Moreover,
  \[
  \Norm1{\vphi_{F,y}}^2=\int_{F^\circledast}\overline{\mathcal F_{F^\perp}(\vphi^x)(y)}\cdot\mathcal F_{F^\perp}(\vphi^x)(y)\,dx
    =\int_{F^*}e^{-2\pi i\Rscp0yx}(\vphi*_{F^\perp}\vphi^{p_F(x)*})(x)\,dx
  \]
  Since the Fourier transform is continuous $\Lp[^1]0X\to\Ct[_0]0X\,$, we need to see that
  \[
  1_{F^*}(x)\cdot(\vphi*_{F^\perp}\vphi^{p_F(x)*})(x)=1_{F^*}(x)\cdot\int_{F^\perp}\vphi(x-w)\overline{\vphi(p_F(x)-w)}\,dw\ ,
  \]
  viewed as an $\mathrm L^1$ function in $x\,$, depends continuously on $(F,y)\,$. This follows from Lebesgue's theorem once we have point-wise continuous dependence, which
  is ensured by Lemma \ref{partial-ft-contin} below.
\end{pf}

\begin{lem}\label{partial-ft-contin}
  Let $\vphi,\psi\in\Cc0X\,$. Define 
  \[
  \chi(F,u,v)=\Bracks0{\vphi*_{F^\perp}\psi^u}(u+v)\mathfa(F,u,v)\in\mathcal W_\Omega|Y_j\ .
  \]
  Then $\chi\in\Cc0{\mathcal W_\Omega|Y_j}\,$.
\end{lem}

\begin{pf}
  Clearly, $\chi$ has compact support, and we need to prove continuity. Let $m=n-n_{d-j}$ denote the common dimension of $F^\perp\,$, $F\in P_j\,$.
  Let $(F_k,u_k,v_k)$ tend to $(F,u,v)\,$, and set
  \[
  \phi_k(w)=1_{F_k^\perp}(w)\cdot\vphi(u_k+v_k-w)\psi(u_k+w)\mathfa k\in\nats\,,\,w\in X\ .
  \]
  Then $\chi(F_k,u_k,v_k)=\int\phi_k\,d\mathcal H^m$ and $\phi_k(w)\to 1_{F^\perp}(w)\cdot\vphi(u+v-w)\psi(u+w)$ \fa $w\in X\,$. (Here, $\mathcal H^m$ denotes 
  $m$-dimensional Hausdorff measure.) There exist $r>0\,$, $C>0$ \scth $\Abs0{\phi_k}\sle C\cdot 1_{F_k^\perp\cap\ball_r}\,$. Note that 
  $\mathcal H^m(F_k^\perp\cap\ball_r)$ is independent of $k\,$, since the intersections are just the $m$-dimensional balls of radius $r$ in $F_k^\perp\,$, centred at the 
  origin. Hence, Pratt's lemma, \cite[th.~1.3.4]{evans-gariepy}, implies that $\chi(F_k,u_k,v_k)\to\chi(F,u,v)\,$.
\end{pf}

\begin{thm}\label{subquot-trivfield}
  The representation $\sigma_j=(L^{F,y})_{(F,y)\in\Sigma_j}$ exhibits $\Cst0{\mathcal W_\Omega|Y_j}$ as isomorphic to the field of elementary C$^*$-algebras
  $\Cpt0{\mathcal E_j}$ associated to $\mathcal E_j\,$. Moreover, this field is trivial.
\end{thm}

\begin{pf}
  Let $\vphi\in\Cc0{\mathcal W_\Omega|Y_j}\,$, and $E\in P_j\,$. By Proposition \ref{norm-loc-triv-ex}, we may choose a normalised local trivialisation $(\psi_U,U)$ at $E\,$.
  Proposition \ref{equiv-rep} shows that $\psi_F^*\Parens1{L^{F,y}(\vphi)}$ depends continuously on $(F,y)\in\Sigma_j|U\,$.
  In particular, $(F,y)\mapsto\Norm0{L_\Omega^{F,y}(\vphi)}$ is continuous. By Proposition \ref{groupoid-loc-triv}, the image of
  $L_\Omega^{F,y}$ on $\Cst0{\mathcal W_\Omega|\pi^{-1}(U)}$ is $\Cst0{\mathcal W_{E^\circledast}|E^\circledast}\cong\Cpt1{\Lp[^2]0{E^\circledast}}\,$.

  By a partition of unity argument, $A=\sigma_j(I_{j+1}/I_j)$ is a locally trivial, continuous field of elementary C$^*$-algebras. It is clear that $\sigma_j$ is injective on
  $\Cst0{\mathcal W_\Omega|Y_j}$, so it sets up an isomorphism with $A\,$.

  To see that the C$^*$-algebra of the continuous field $\mathcal E_j$ is contained in $A\,$, it suffices to
  see that $\vtheta_{\vphi,\psi}:(F,y)\mapsto\psi_{F,y}^{\vphantom*}\vphi_{F,y}^*$ lies in $A$ \fa $\vphi,\psi\in\Cc0X\,$. Let $s_F=2p_F-1$ and
  \[
  \chi(F,u,v)=\Bracks0{(\bar\vphi\circ s_F)*_{F^\perp}\psi^u}(u+v)\mathfa(F,u,v)\in Y_j\ .
  \]
  Then $\chi\in\Cc0{\mathcal W_\Omega|Y_j}$ by Lemma \ref{partial-ft-contin}. Now,
  \begin{align*}
    L^{F,y}(\chi)h(u)&=\int_{F^\circledast}\mathcal F_{F^\perp}\Parens1{\Bracks0{(\bar\vphi\circ s_F)*_{F^\perp}\psi^u}^v}(y)h(v)\,dv\\
    &=\mathcal F_{F^\perp}(\psi^u)(y)\cdot\int_{F^\circledast}\overline{\mathcal F_{F^\perp}(\vphi^v)(y)}\cdot h(v)\,dv=\vtheta_{\vphi,\psi}h(u)\ .
  \end{align*}
  This implies $\Cpt0{\mathcal E_j}\subset A\,$, and since the former separates points, equality, by \cite[lem.~10.5.3]{dixmier}.

  The triviality of the field $\mathcal E_j$ for $j=d$ is clear, since $\Sigma_d\approx X$ is contractible. For $j<d\,$, it follows from \cite[lem.~10.8.7]{dixmier} since its
  fibre $\Lp[^2]0{F^\circledast}$ is separable, and its base $\Sigma_j$ is finite-dimensional by Lemma \ref{face-space-fd} below.
\end{pf}

\begin{lem}\label{face-space-fd}
  For $0\sle j\sle d\,$, the map $P_j\to\Gr_{n_{j-d}}(X):F\mapsto\Span0F$ is a topological embedding into the Grassmannian of $n_{j-d}$-planes. Consequently, the spaces
  $P_j$ and $\Sigma_j$ are finite-dimensional.
\end{lem}

\begin{pf}
  The map is continuous by Proposition \ref{lower-dim-closed}, and injective since $F=\Omega^*\cap\Span0F\,$. Thus, it is topological, seeing that $P_j$ is compact.
  The image of $P_j$ has dimension $\sle n_{d-j}\cdot\Parens1{n-n_{d-j}}\,$, by \cite[ch.~III, \S~1, th.~III.1]{hurewicz-wallmann}.  Moreover, dimension is invariant
  under homeomorphisms, cf.~\cite[ch.~III, \S~1, rem.~A)]{hurewicz-wallmann}. The finite-dimensionality of $P_j$ follows, and \cite[ch.~III, \S~4, th.~III.4]{hurewicz-wallmann}
  entails that of $\Sigma_j\subset P_j\times X\,$.
\end{pf}

\begin{rem}
  Needless to say, our proof of Theorem \ref{subquot-trivfield} follows the proof of the corresponding results in \cite[th.~4.7, th.~6.4.]{muhly-renault} for polyhedral and symmetric
  cones quite closely; the main new ingredient being the application of some convex analysis to the construction of local trivialisations.
\end{rem}

\begin{cor}
  For $0\sle j<d\,$, the C$^*$-algebra $\Cst0{\mathcal W_\Omega|Y_j}$ is stable.
\end{cor}

\begin{pf}
  Indeed, the trivial field $\mathcal E_j$ has separable infinite-dimensional fibre for $j<d\,$.
\end{pf}

\subsection{Analytical Index Formula}
For any Hilbert C$^*$-module $\mathcal E\,$, let $Q(\mathcal E)=\Lct0{\mathcal E}/\Cpt0{\mathcal E}$ denote its Calkin algebra. Let 
$\tau_j:\Cst0{\mathcal W_\Omega|Y_j}\to Q(\mathcal E_j)$ be the Busby invariant of the extension from section \ref{jth-wh-ext}. We call the this the
  \emph{$j$th Wiener--Hopf extension}. 

  If $\vrho_j:\Cst0{\mathcal W_\Omega}\to I_{j+1}/I_{j-1}$ is a completely positive contractive section of $\sigma_j\,$, then
  $\tau_j=q_{j-1}\circ\sigma_{j-1}\circ\vrho_j$ where $q_{j-1}:\Lct0{\mathcal E_{j-1}}\to Q(\mathcal E_{j-1})$ is the canonical projection onto the Calkin algebra of the
  Hilbert $\Ct[_0]0{\Sigma_j}$-module $\mathcal E_{j-1}\,$, and $\sigma_{j-1}:I_{j+1}/I_{j-1}\to\Lct0{\mathcal E_{j-1}}$ is the strict extension of
  $\sigma_{j-1}:\Cst0{\mathcal W_\Omega|Y_{j-1}}\to\Cpt0{\mathcal E_{j-1}}\,$.

  Moreover, by naturality of connecting homomorphisms, $\sigma_{j-1}^*\partial_j=\tau_{j*}\partial$ where $\partial$ is the connecting homomorphism of
  \[
  \begin{CD}
    0@>>>{\Cpt0{\mathcal E_{j-1}}}@>>>{\Lct0{\mathcal E_{j-1}}}@>{q_j}>>{Q(\mathcal E_{j-1})}@>>>0\ .\end{CD}
  \]

  \par\noindent
  We call an element $a\in\Cst0{\mathcal W_\Omega}$ \emph{$j$-Fredholm} if it represents an invertible in the unitisation of the quotient
  $\Cst0{\mathcal W_\Omega}/I_j\,$. Equivalently, $ab\equiv ba\equiv 1\pmod{I_j}$ \fs $b\in\Cst0{\mathcal W_\Omega}\,$. More generally, any
  $a\in\Cst0{\mathcal W_\Omega}\otimes\cplxs^{N\times N}$ which is invertible modulo $I_j\otimes\cplxs^{N\times N}$ shall be called a $j$-\emph{Fredholm matrix}.

\begin{prop}\label{fred-family}
  If $a\in\Cst0{\mathcal W_\Omega}$ is $j$-Fredholm, then $\sigma_{j-1}(a)=\Parens1{L^{F,y}(a)}_{(F,y)\in\Sigma_{j-1}}$ is a continuous family of Fredholm
  operators. The corresponding statement about matrices is also valid.
\end{prop}

\par\noindent
  To that end, we observe the following naturality of the representations $L^{F,y}\,$.

\begin{lem}
  Let $F\in P\,$, and let $P_F$ be the set of faces of $F\,$. Then we may define a $*$-homomorphism $r_F:\Cst0{\mathcal W_\Omega}\to\Cst0{\mathcal W_{F^\circledast}}$ by
  \[
  r_F(\vphi)(E,u,v)=\int_{F^\perp}\vphi(E,u,y+v)\,dy\mathfa(E,u,v)\in\mathcal W_{F^\circledast}\,,\,\vphi\in\Cc0{\mathcal W_\Omega}\ .
  \]
  Moreover, we have
  \[
  L^{E,v}_{F^\circledast}\circ r_F=L^{E,v}_\Omega\mathfa E\in P_F\,,\,v\in E^\perp\cap\Span0F\ .
  \]
\end{lem}

\begin{pf}
  Observe $E^*=\Span0F\cap E^*\oplus F^\perp\,$, since $F^\perp\subset E^\perp$ \fa $E\in P_F\,$. For $\vphi\in\Cc0{\mathcal W_\Omega}\,$, we compute
  \begin{align*}
    L^{E,v}_{F^\circledast}r_F(\vphi)h(u)&=\int_{\Span0F\cap E^*-u}\int_{F^\perp}\vphi(E,u,y+w)e^{-i\Rscp0{w+y}v} h\Parens1{u+p_E(w+y)}\,dy\,dw\\
    &=\int_{E^*-u}\vphi(E,u,w)e^{-i\Rscp0wv}h\Parens1{u+p_E(w)}\,dw=L^{E,v}_\Omega(\vphi)h(u)
  \end{align*}
  \fa $E\in P_F\,$, $v\in\Span0F\cap E^\perp\,$, $h\in\Lp[^2]0{E^\circledast}\,$, and $u\in E^\circledast\,$, since in the integral, $y$ is perpendicular to $v,w\,$. This proves 
  the second equality. Choosing $E=F\,$, $v=0\,$, $L^{F,0}_{F^\circledast}$ is an isomorphism onto its image, so $r_F$ is bounded, and an involutory algebra
  homomorphism. This proves the lemma.
\end{pf}

{\def\Elproofname{\protect{PROOF of Proposition \ref{fred-family}.}}
\begin{pf}
  The statement about Fredholm matrices follows along the same lines as the first assertion, so for the sake of simplicity, we restrict ourselves to $N=1\,$.
  Then the continuous dependence is clear from Proposition \ref{cont-field}.

  Let $b\in\Cst0{\mathcal W_\Omega}\,$, $ab\equiv ba\equiv 1\pmod{I_j}\,$. Take $(F,y)\in\Sigma_{j-1}\,$, and let
  $E\in P_F\,$, $E\neq F\,$. Then $E\in P\,$, and hence $\dim E\sle n_{d-j}\,$. Thus,
  \[
  I_j\subset\ker L^{E,v}_\Omega\mathfa v\in\Span0F\cap E^\perp\ .
  \]
  This implies
  \[
  1=L^{E,v}_\Omega(ab)=L^{E,v}_{F^\circledast}r_F(ab)\mathfa v\in\Span0F\cap E^\perp\ .
  \]
  Since $E$ was arbitrary,
  \[
  r_F(ab)-1\in\bigcap\nolimits_{E\in P_F\setminus\{F\}\,,\,v\in\Span0F\cap E^\perp}\ker L^{E,v}_{F^\circledast}
  =\Parens1{L^{F,0}_{F^\circledast}}^{-1}\Parens1{\Cpt0{\protect{\Lp[^2]0{F^\circledast}}}}\ ,
  \]
  by the composition series for $\Cst0{\mathcal W_{F^\circledast}}\,$. We conclude
  \[
  L^{F,0}_\Omega(a)L^{F,0}_\Omega(b)-1=L^{F,0}_\Omega(ab)-1=L^{F,0}_{F^\circledast}\Parens1{r_F(ab)-1}\in\Cpt0{\protect{\Lp[^2]0{F^\circledast}}}\ .
  \]
  If we denote by $e^{-iy^*}$ the bounded continuous function $\mathcal W_\Omega\to\cplxs:(E,u,v)\mapsto e^{-i\Rscp0yv}\,$, then
  $L^{F,y}_\Omega(\vphi)=L^{F,0}_\Omega(e^{-iy^*}\cdot\vphi)\,$. Thus, the above entails
  \[
  L^{F,y}_\Omega(a)L^{F,y}_\Omega(b)-1=L^{F,0}_\Omega(e^{-iy^*}\cdot ab)-1=L^{F,0}_{F^\circledast}(r_F(e^{-iy^*}\cdot ab)-1)
  \in\Cpt0{\protect{\Lp[^2]0{F^\circledast}}}\ .
  \]
  Similarly, $L^{F,y}_\Omega(b)L^{F,y}_\Omega(a)-1$ is compact. Hence, $L^{F,y}_\Omega(a)$ is Fredholm.
\end{pf}}

\par\noindent
  Recall that $[f]\in K^1_c(\Sigma_j)$ is given by a continuous map $f:\Sigma_j\to\mathrm U(N)$ for some $N\in\nats\,$, such that $(f_{k\ell})=(\delta_{k\ell})$
  outside some compact set. Fix some completely positive cross section $\vrho_j:\Cpt0{\mathcal E_j}\to I_{j+1}$ of $\sigma_j\,$. We claim that
  \[
  \vrho_j(f)=1_N+\Parens1{\vrho_j\Parens1{f_{k\ell}-\delta_{k\ell}}}_{1\sle k,\ell\sle N}
  \]
  is a $j$-Fredholm matrix. (Here, we identify $f-1_N\in\Ct[_0]0{\Sigma_j}\otimes\cplxs^{N\times N}$ with its preimage in $\Cpt0{\mathcal E_j\otimes\cplxs^N}\,$.)
  Indeed, denoting the unital extension of $\sigma_j\otimes\id_{\cplxs^{N\times N}}$ to the unitisation of $I_{j+1}\otimes\cplxs^{N\times N}$ by $\sigma_j\,$,
  \[
  \sigma_j\Parens1{1+\vrho_j(f-1)}=1+\sigma_j\vrho_j(f-1)=f\ ,
  \]
  which is invertible in the unitisation of $\Cpt0{\mathcal E_j\otimes\cplxs^N}\,$. Since $\ker\sigma_j=I_j\,$, this means that $\vrho_j(f)$ is a $j$-Fredholm matrix.

  It is therefore natural to ask whether the map $\partial_j$ can be interpreted as the Atiyah--J\"anich family index of the family
  $\sigma_{j-1}\vrho_j(f)=1+\sigma_{j-1}\vrho_j(f-1)$ of Fredholm operators. First, we need to see that such a family index is well-defined.

\begin{prop}\label{fred-fam-van-infty}
  Let $[f]\in K^1_c(\Sigma_j)\,$. Then $\sigma_{j-1}\vrho_j(f)$ is trivial at infinity, i.e.~there exists a compact $L\subset\Sigma_{j-1}$ \scth
  \[
  L^{F,y}\vrho_j(f_{k\ell}-\delta_{k\ell})\in\Cpt0{\protect{\Lp[^2]0{F^\circledast}}}\mathfa(F,y)\in\Sigma_{j-1}\setminus L\ .
  \]
\end{prop}

\par\noindent
  We first make the following observation. Let $\mathcal P$ denote the graph of the order $\supset$ of the face lattice $P\,$, and
  \[
  \mathcal P_j=\mathcal P\cap(P_{j-1}\times P_j)=\Set1{(E,F)\in P_{j-1}\times P_j}{E\supset F}\ .
  \]
  Moreover, denote its projections by $\begin{CD}P_{j-1}@<\xi<<{\mathcal P_j}@>\eta>>&P_j\end{CD}\,$.

\begin{lem}\label{sigmas-cpt}
  The relation $\mathcal P$ is closed. Thus, $\mathcal P_j$ is a compact subspace of $P_{j-1}\times P_j\,$. The projections $\xi$ and $\eta$ are continuous, closed, and proper.
\end{lem}

\begin{pf}
  Let $(E_k,F_k)\in\mathcal P\,$, $(E_k,F_k)\to(E,F)\in P\times P\,$. If $e\in E\,$, then $e=\limk_ke_k$ \fs $e_k\in E_k\subset F_k\,$. Hence, $e\in\limk_kF_k=F\,$. Therefore,
  $\mathcal P$ is closed. The continuity of $\xi$ and $\eta$ is clear. The closedness and properness follow from the compactness of $\mathcal P_j\,$.
\end{pf}

{\def\Elproofname{\protect{PROOF of Proposition \ref{fred-fam-van-infty}.}}
\begin{pf}
  Let $[f]\in K^1_c(\Sigma_j)$ where $f:\Sigma_j\to\mathrm U(N)$ \fs $N\in\nats$ and $f=1_N$ on $\Sigma_j\setminus K$ where $K$ is compact. Since $\Sigma_j$ is a vector
  bundle over $P_j\,$, we may consider $\eta^*K\subset\eta^*\Sigma_j\,$. Due to the properness of $\eta\,$, this set is compact. The projection
  \[
  \eta^*\Sigma_j\to\xi^*\Sigma_{j-1}:(E,F,y)\mapsto(E,p_{E^\perp}(y),F)
  \]
  is continuous, so we obtain a compact subset of $\xi^*\Sigma_{j-1}$ which is necessarily of the form $\xi^*L$ for some compact $L\subset\Sigma_{j-1}\,$. Explicitly, $L$ may
  be written down as follows,
  \[
  L=\Set1{(E,v)\in\Sigma_{j-1}}{\exists\,F\in\eta(\xi^{-1}(E))\,,\,u\in F^\perp\cap\Span0E\,:\,(F,u+v)\in K}\ .
  \]

  Fix $(E,v)\in\Sigma_{j-1}\setminus L\,$.  Let $H\in P_E\,$, $H\neq E\,$, and $w\in H^\perp\cap\Span0E\,$.
  We have $F\not\subset E$ for every $F\in P_j\,$, $u\in F^\perp\cap\Span0E$ \scth $(F,u+v)\in K\,$. On the other hand, $H\subset E\,$, \sth $(H,v+w)\not\in K\,$. Hence,
  \begin{align*}
    L_{E^\circledast}^{H,v}\Parens1{L_{E^\circledast}^{E,0}}^{-1}L_\Omega^{E,y}\vrho_j(f_{k\ell}-\delta_{k\ell})
    &=L_{E^\circledast}^{H,v}r_E\Parens1{e^{-iy^*}\cdot\vrho_j(f_{k\ell}-\delta_{k\ell})}\\
    &=L_\Omega^{H,v+y}(\vrho_j(f_{k\ell}-\delta_{k\ell}))=\Parens1{f_{k\ell}-\delta_{k\ell}}(H,v+y)=0\ .
  \end{align*}
  Thus,
  \[
  \Parens1{L_{E^\circledast}^{E,0}}^{-1}L_\Omega^{E,y}\vrho_j(f_{k\ell}-\delta_{k\ell})\in\bigcap_{(H,v)}\ker L_{E^\circledast}^{H,v}
  =\Parens1{L_{E^\circledast}^{E,0}}^{-1}\Parens1{\Cpt0{\protect{\Lp[^2]0{E^\circledast}}}}
  \]
  \fa $(E,y)\in\Sigma_{j-1}\setminus L\,$, which proves our assertion.
\end{pf}}

\par\noindent
 Proposition \ref{fred-fam-van-infty} enables us to define the Atiyah--J\"anich family index of the continuous family $\sigma_{j-1}\vrho_j(f)$ of Fredholm operators,
  where $f\in K^1_c(\Sigma_j)\,$, by the following standard device. Consider a filtration $X_0\subset X_1^\circ\subset X_1\subset\dotsm\subset\Sigma_{j-1}$ by compact
  sets whose interiors $X_k^\circ$ are non-void and whose union is $\Sigma_{j-1}\,$. For each $k\in\nats\,$, the index
  \[
  \Index\sigma_{j-1}\vrho_j(f)\bigm|X_k\in K^0(X_k)=K^0_c(X_k)
  \]
  is well-defined, cf.~\cite[p.~158]{atiyah-ktheory}, \cite[p.~138]{jaenich-index}.

  Let $\xi_k=\Index\sigma_{j-1}\vrho_j\bigm|X_k^\circ\in K_c^0(X_k^\circ)\,$, and denote by $j_k:K^0_c(X_k^\circ)\to K^0_c(X_{k+1}^\circ)$ the respective
  wrong way maps (i.e.~extension by zero). Then $j_k(\xi_k)=\xi_{k+1}$ for $k$ large enough, since outside some $X_k^\circ\,$, $\sigma_{j-1}\vrho_j(f)$ is trivial.
  If we write $T=\sigma_{j-1}\vrho_j(f)\,$, this means $T_{F,y}=1_N$ for $(F,y)\not\in X_k\,$, possibly replacing $T$ by a homotopic family
  (the set of compact operators is convex). But then $T_{F,y}(V)=V$ for any $V$ of finite codimension.
  By construction of the family index (loc.~cit.), this shows that the restriction of $\xi_\ell$ to $\Sigma_{j-1}\setminus X_k^\circ$ vanishes for $\ell>k\,$. By naturality of the
  index, \cite[p.~159]{atiyah-ktheory}, \cite[lem.~6]{jaenich-index}, the restriction of $\xi_\ell$ to $X_k^\circ$ is $\xi_k\,$.
  Thus, we indeed have $j_k(\xi_k)=\xi_{k+1}\,$. Since $K^0_c(\Sigma_{j-1})=\varinjlim\nolimits_k K^0_c(X_k^\circ)\,$, by \cite[ch.~II, prop.~4.21]{karoubi-ktheory}, we find
  that there exists a uniquely determined $\xi\in K^0_c(\Sigma_{j-1})$ such that its restriction to $X_k^\circ$ is $\xi_k\,$. We denote the class $\xi$ by 
  $\Index_{\Sigma_{j-1}}\sigma_{j-1}\vrho_j(f)\,$. 

\begin{thm}
  For $[f]\in K^1_c(\Sigma_j)\,$, we have
  \[
  \partial_j[f]=\Index_{\Sigma_{j-1}}\sigma_{j-1}\vrho_j(f)\ ,
  \]
  for any choice of completely positive contractive section $\vrho_j:\Cpt0{\mathcal E_j}\to\Cst0{\mathcal W_\Omega|U_{j+1}}$ for $\sigma_j\,$.
\end{thm}

\begin{pf}
  By naturality of connecting maps, it suffices to establish the fact that the connecting map for the extension
  \[
  \begin{CD}0@>>>{A\otimes\mathbb K}@>>>{\Mult0{A\otimes\mathbb K}}@>q>>Q(A\otimes\mathbb K)@>>>0\ ,\end{CD}
  \]
  where $A=\Ct0Z$ for some compact space $Z\,$, is given by the Atiyah--J\"anich index. This follows exactly as for $Z$ a point. Indeed, let 
  $[u]\in K_1(Q(A\otimes\mathbb K))\,$, \sth 
  \[
  u^*u\equiv uu^*\equiv1\pmod{A\otimes\mathbb K}\ .
  \]

  By \cite[prop.~1.5, prop.~1.7]{mingo-fred}, there exists a partial isometry $v\in M(A\otimes\mathbb K)$ \scth we have $u-v\in A\otimes\mathbb K\,$, and 
  $1-vv^*$ and $1-v^*v$ have finitely generated range. So the ranges are contained in the range of the standard projection 
  $p_N:A\otimes\ell^2\to A\otimes\cplxs^N$ for $N\gg0\,$. Then
  \[
  \Index [u]=[1-v^*v]-[1-vv^*]=[wp_Nw^{-1}]-[p_N]=\partial[u]\smathtxt{where}w=\begin{Matrix}0v&1-vv^*\\1-v^*v&v^*\end{Matrix}\ ,
  \]
  which proves the theorem.
\end{pf}

\begin{rem}
  The above deduction of the analytic expression of the index maps $\partial_j$ owes much to the exposition of \cite{upmeier_toeplitz} of the index maps for
  Toeplitz operators; the main differences again being the reconstruction of the Jordan algebraic computations performed there in terms of the convex geometry of the
  cone, and of course the groupoid framework for the C$^*$-algebras involved. Let us remark that our proof of the topological index formula in \cite{alldridge-johansen-wh2} 
  uses methods completely different from Upmeier's, and in particular, contains as a special case an independent proof of the index formula from \cite{upmeier-wh} for 
  symmetric cones.
\end{rem}

\end{document}